\documentclass[11pt]{amsart}
\usepackage{amsmath,graphicx}
\usepackage[all]{xy}

\newtheorem{thm}{Theorem}[section]
\newtheorem{lem}[thm]{Lemma}
\newtheorem{prop}[thm]{Proposition}
\newtheorem{cor}[thm]{Corollary}

\theoremstyle{definition}

\theoremstyle{remark}

\theoremstyle{plain}

\newcommand{\Z}{{\mathbb{Z}}}

\newcommand{\R}{{\mathbb{R}}}
\newcommand{\D}{\mathbb D}

\newcommand{\FF}{{\mathcal{F}}}

\newcommand{\DD}{{\mathcal{D}}}
\newcommand{\RR}{\mathcal{R}}
\newcommand{\eps}{\varepsilon}

\newcommand{\area}{{\mathrm{Area}}}
\newcommand{\dist}{{\mathrm{dist}}}

\newcommand{\II}{\mathcal{I}}
\newcommand{\JJ}{\mathcal{J}}

\newcommand{\EE}{\mathcal{E}}
\newcommand{\NN}{\mathcal{N}}
\newcommand{\MM}{\mathcal{M}}
\newcommand{\MMM}{\widetilde{\MM}}

\numberwithin{equation}{section}

\begin{document}

\title[]{On the correlations of directions in the Euclidean plane}

\author[]{Florin P. Boca and Alexandru Zaharescu}

\address{ Department of Mathematics, University of Illinois, Urbana IL 61801, USA}

\address{Institute of Mathematics "Simion Stoilow" of the Romanian
Academy, P.O. Box 1-764, RO-014700 Bucharest, Romania}

\address{E-mail: fboca@math.uiuc.edu; zaharesc@math.uiuc.edu}

\date{April 5, 2004.}

\begin{abstract}
Let $\RR^{(\nu)}_{(x,y),Q}$ denote the repartition of the
$\nu$-point correlation measure of the finite set of directions
$P_{(x,y)}P$, where $P_{(x,y)}$ is the fixed point $(x,y)\in
[0,1)^2$ and $P$ is an integer lattice point in the square
$[-Q,Q]^2$. We show that the average of the pair correlation
repartition $\RR^{(2)}_{(x,y),Q}$ over $(x,y)$ in a fixed disc
$\D_0$ converges as $Q\rightarrow \infty$. More precisely we
prove, for every $\lambda \in \R_+$ and each
$0<\delta<\frac{1}{10}$, the estimate
\begin{equation*}
\frac{1}{\mathrm{Area} (\D_0)} \iint\limits_{\D_0}
\RR^{(2)}_{(x,y),Q}(\lambda)\, dx\, dy=\frac{2\pi
\lambda}{3}+O_{\D_0,\lambda,\delta} (Q^{-\frac{1}{10}+\delta})
\qquad \mbox{\rm as $Q\rightarrow \infty$.}
\end{equation*}
We also prove that for each individual point $(x,y)\in [0,1)^2$,
the $6$-level correlation $\RR^{(6)}_{(x,y),Q}(\lambda)$ diverges
at any point $\lambda \in \R^5_+$ as $Q\rightarrow \infty$, and
give an explicit lower bound for the rate of divergence.
\end{abstract}

\maketitle

\section{Introduction}
In many problems one is led to consider in the Euclidean plane
lines joining a fixed point $P_0$ (which is not necessarily an
integer lattice point) with a set of integer lattice points. A
natural way of measuring the distribution of directions $P_0 P$,
$P\in \Z^2$, is via correlations and consecutive spacings. When
the fixed point is the origin, the problem is related to the
distribution of Farey fractions with multiplicities, each fraction
$\frac{a}{q}$ in $\FF_Q$ being counted $\big[ \frac{Q}{q}\big]$
times. The consecutive $h$-level spacing measures of customary
Farey fractions were computed for $h=1$ in \cite{Hall} and for
$h\geq 2$ in \cite{ABCZ}. Limiting correlations of Farey fractions
were shown to exist and computed recently in \cite{BZ2}.

When the fixed point is not an integer lattice point, the problem
of existence of limiting correlations/consecutive spacings is
considerably more difficult. It is therefore natural to try to
prove first some averaging results, letting the fixed point to
vary in a given region. In the first part of this paper we derive
such a result for the limiting pair correlation measure. The
limiting average pair correlation function is constant, as in the
Poisson case. What is striking however is that this constant is
not $1$, as in the Poisson case, but $\frac{\pi}{3}$.

We now give a mathematical formulation of the problem. For each
$Q\geq 1$, let $\Box_Q$ denote the set of integer lattice points
in the square $[-Q,Q]^2$, and set $N=N_Q=\# \Box_Q=(2Q+1)^2$. For
each point $P_{(x,y)}=(x,y)$, we consider the finite sequences
$\big( \theta_P(x,y)\big)_{P\in \Box_Q}$, $Q$ large integer, of
angles between the line $P_{(x,y)}P$ and the horizontal direction.
The pair correlation of this finite sequence is defined as
\begin{equation*}
\RR^{(2)}_{(x,y),Q}(\lambda) =\frac{\#\left\{ (P,P^\prime)\in
\Box_Q^2 \, :\, P\neq P^\prime,\ \frac{N}{2\pi} \, \vert
\theta_{P,P^\prime}(x,y)\vert \leq \lambda \right\}}{N}\ ,\qquad
\lambda \in \R_+,
\end{equation*}
where $\theta_{P,P^\prime}(x,y)$ denotes the measure of the angle
$\angle PP_{(x,y)}P^\prime$.

Throughout the paper we shall consider a fixed disc $\D_0$ of
center $(x_0,y_0)\in [0,1)^2$ and radius $r_0$. We are interested
in the asymptotic behavior of the average
\begin{equation*}
R^{(2)}_{\D_0,Q} (\lambda)=\frac{1}{\pi r_0^2} \iint\limits_{\D_0}
\RR^{(2)}_{(x,y),Q}(\lambda)\, dx\, dy
\end{equation*}
of $\RR^{(2)}_{(x,y),Q}(\lambda)$ over $\D_0$, for fixed $\lambda
>0$ and $Q\rightarrow \infty$.

\medskip

The first three sections are concerned with the proof of the
following result.
\begin{thm}\label{T1.1}
For every $\lambda >0$ and $\delta >0$
\begin{equation}\label{T1}
R_{\D_0,Q}^{(2)}(\lambda)=\frac{2\pi
\lambda}{3}+O_{\D_0,\lambda,\delta} (
Q^{-\frac{1}{10}+\delta})\qquad \mbox{as $\ Q\rightarrow \infty$.}
\end{equation}
\end{thm}

\medskip

If one replaces $\D_0$ by a vertical or horizontal segment of
length one, an identical asymptotic formula as in \eqref{T1} turns
out to be true. This can be proved by similar techniques as in
this paper or using Erd\" os-Tur\' an type discrepancy estimates,
and suggests that \eqref{T1} may be true regardless of the shape
of the range of the fixed point.

The behavior of higher level correlations appears to be different.
In the last section we prove that the $6$-level correlations
diverge for every \underline{individual} fixed point. When $\nu
\geq 2$, the repartition of the $\nu$-level correlation measure of
the finite sequence $\big( \theta_P (x,y)\big)_{P\in \Box_Q}$ is
defined for each vector $\lambda =(\lambda_1,\dots,\lambda_{\nu
-1})\in \R_+^{\nu -1}$ by
\begin{equation}\label{1.2}
\mbox{\small $\displaystyle
\RR^{(\nu)}_{(x,y),Q}(\lambda)=\frac{\# \left\{
(P_1,\dots,P_\nu)\in \Box_Q^\nu\, :\ \mbox{\rm $P_i$ distinct},\
\vert \theta_{P_i,P_{i+1}}(x,y)\vert \leq \frac{2\pi \lambda_i}{N}
\, ,\ 1\leq i\leq \nu-1\right\}}{N}\ .$}
\end{equation}

For randomly chosen directions one would expect to obtain the
Poissonian limit
\begin{equation}\label{1.3}
\lim\limits_{Q\rightarrow \infty} \RR^{(\nu)}_{(x,y),Q}(\lambda)
=\mathrm{Vol} \ \prod\limits_{i=1}^{\nu-1} \
[-\lambda_i,\lambda_i] = 2^{\nu-1} \lambda_1 \dots \lambda_{\nu-1}
.
\end{equation}
It turns out however that \eqref{1.3} fails in this situation.
More precisely, we will show that if $\nu\geq 6$, then for every
point $(x,y)\in [0,1)^2$ and for every
$(\lambda_1,\dots,\lambda_{\nu-1})\in \R_+^{\nu-1}$,
$\lim\limits_{Q\rightarrow \infty}
\RR^{(2)}_{(x,y),Q}(\lambda)=\infty$. This is a consequence of

\begin{thm}\label{T1.2}
For every $(x,y)\in [0,1)^2$, every $\lambda
=(\lambda_1,\dots,\lambda_5)\in \R_+^5$, and every $\delta>0$, for
$Q$ large enough in terms of $x,y,\lambda$ and $\delta$
\begin{equation}\label{IV.28}
\RR^{(6)}_{(x,y),Q}(\lambda)>Q^{\frac{1}{4}-\delta} .
\end{equation}
\end{thm}

As in Theorem \ref{T1.2} one can prove
\begin{cor}\label{C1.3}
The $6$-level correlations of angles of directions $P_{(x,y)}P$,
where $P$ is a lattice point inside an expanding region $Q\Omega$,
diverges as $Q\rightarrow \infty$ whenever $\Omega$ is a convex
domain in $\R^2$ which contains the origin.
\end{cor}

The phenomenon is similar to the one encountered in the problem of
the distribution of fractional parts of polynomials. There, one
can handle the pair correlation problem generically (see
\cite{RS},\cite{BZ}). Moreover, in the case of the sequence $n^2
\alpha \pmod{1}$ one is able to solve the problem for all
$m$-level correlations for a large class of irrational numbers
$\alpha$ (see \cite{RSZ},\cite{Z}). However, as shown in
\cite{RSZ}, there are irrational numbers $\alpha$ for which the
$5$-level correlation of fractional parts of $n^2 \alpha$, $1\leq
n\leq N$, diverges to infinity as $N\rightarrow \infty$. This
occurs as a result of the presence of large clusters of such
fractional parts. In the case of Theorem \ref{T1.2} above, large
clusters of elements of the given sequence are responsible, too,
for the divergence of the $6$-level correlations, and hence of any
other higher level correlations.

\bigskip

\section{A first approximation for $R^{(2)}_Q (\lambda)$}

For obvious practical reasons, we try to replace from the
beginning $\theta_{P,P^\prime} (x,y)$ by one of its trigonometric
functions in the definition of $\RR^{(2)}_{(x,y),Q} (\lambda)$.
Suppose that two distinct points $P=(q,a)$, $P^\prime =(q^\prime
,a^\prime)\in\Box_Q$, are such that $q,q^\prime \geq 0$ and $\max
\{ a,a^\prime \}> 0>\min \{ a,a^\prime \}$. Then for sufficiently
large $Q$ (depending only on $\lambda$) we have
\begin{equation*}
\min_{(x,y)\in [0,1]^2}\vert \theta_{P,P^\prime}(x,y)\vert \geq
\arcsin \frac{1}{\sqrt{Q^2+1}} >\frac{2\pi \lambda}{N} \ .
\end{equation*}

As a result, we may only consider in the definition of
$\RR^{(2)}_{(x,y),Q}$ points from the same quadrant. Thus if we
set
\begin{equation*}
\tilde{\Box}_Q^2=\Big\{ (P,P^\prime)\in \Box_Q^2 \, :\, P\neq
P^\prime \ \mbox{\rm and}\ P,P^\prime\ \mbox{\rm belong to the
same quadrant}\Big\}
\end{equation*}
and
\begin{equation}\label{2.1}
\beta_{Q,\lambda}=\sin \frac{2\pi \lambda}{N}=\sin \frac{2\pi
\lambda}{(2Q+1)^2} =\frac{\pi \lambda}{2Q^2} +O_\lambda\bigg(
\frac{1}{Q^6}\bigg) \qquad \mbox{\rm as $Q\rightarrow \infty$,}
\end{equation}
then
\begin{equation}\label{2.2}
\begin{split}
\RR^{(2)}_{(x,y),Q} (\lambda) & =\frac{\# \left\{ (P,P^\prime)\in
\tilde{\Box}_Q^2 \, :\, \ \vert \theta_{P,P^\prime}(x,y)\vert \leq
\frac{2\pi \lambda}{N} \right\}}{N}
\\ & =\frac{\# \left\{ (P,P^\prime)\in \tilde{\Box}_Q^2 \,
:\, \vert \sin \theta_{P,P^\prime}(x,y)\vert \leq
\beta_{Q,\lambda}\right\}}{N}\  .
\end{split}
\end{equation}

For $P=(q,a),P^\prime =(q^\prime ,a^\prime)$, $(x,y)\in \R^2$, we
define
\begin{equation*}
L_{P,P^\prime}(x,y)=(a^\prime -y)(q-x)-(a-y)(q^\prime -x)=\left|
\begin{matrix} 1 & q & a \\
1 & q^\prime & a^\prime \\
1 & x & y \end{matrix} \right| .
\end{equation*}
Then
\begin{equation*}
\vert \sin \theta_{P,P^\prime}(x,y)\vert =\frac{2\area \triangle
PP_{(x,y)} P^\prime}{\| P_{(x,y)}P\| \cdot \| P_{(x,y)}P^\prime\|}
=\frac{\vert L_{P,P^\prime}(x,y)\vert}{\| P_{(x,y)}P\| \cdot \|
P_{(x,y)}P^\prime\|} \ .
\end{equation*}
For each $P,P^\prime \in \Box_Q$, consider the weight
\begin{equation*}
w_{P,P^\prime}(Q,\lambda)=\area \Big\{ (x,y)\in \D_0 \, :\, \vert
L_{P,P^\prime}(x,y)\vert \leq \beta_{Q,\lambda} \| P_{(x,y)} P\|
\cdot \| P_{(x,y)}P^\prime \| \Big\} .
\end{equation*}
From \eqref{2.2} we infer that
\begin{equation}\label{2.3}
R^{(2)}_{Q} (\lambda) =\frac{1}{N} \sum\limits_{(P,P^\prime)\in
\tilde{\Box}_Q^2} \hspace{-10pt} w_{P,P^\prime} (Q,\lambda) .
\end{equation}
Denote
\begin{equation*}
\gamma=\gamma_{P,P^\prime}(Q) =\frac{\| OP\|\cdot \| OP^\prime
\|}{Q^2} =\frac{\sqrt{q^2+a^2} \ \sqrt{q^{\prime 2}+a^{\prime
2}}}{Q^2}\ ,
\end{equation*}
and define for every $\mu >0$
\begin{equation}\label{2.4}
\begin{split}
& A_{P,P^\prime}(Q,\mu)=\area \Big\{ (x,y)\in \D_0 \, :\, \vert
L_{P,P^\prime}(x,y)\vert \leq \mu
\gamma_{P,P^\prime}(Q)\Big\} ,\\
& G_Q(\mu) =\frac{1}{Q^2} \sum\limits_{(P,P^\prime) \in
\tilde{\Box}_Q^2} A_{P,P^\prime}(Q,\mu) .
\end{split}
\end{equation}

In the remainder of this section we show that the asymptotic of
$R^{(2)}_Q(\lambda)$ as $Q\rightarrow \infty$ is closely related
to that of $G_Q \big( \frac{\pi\lambda}{2}\big)$.

For fixed $P,P^\prime$, denote by $\theta$ the angle between the
line $\ell$ determined by $P$ and $P^\prime$ and the horizontal
direction. Consider also the lines $\ell_\pm$, parallel to $\ell$
and such that $\dist (\ell,\ell_\pm)=\frac{\mu \gamma \cos
\theta}{\vert q^\prime -q\vert}$. The equation of $\ell$ is given
by
\begin{equation*}
\mbox{\rm ($\ell$)} \qquad \qquad \qquad L_{P,P^\prime}(x,y)=0,
\end{equation*}
while the equation of $\ell_\pm$ is given by
\begin{equation*}
\mbox{\rm ($\ell_\pm$)} \qquad \qquad L_{P,P^\prime}(x,y)=\pm \mu
\gamma  .
\end{equation*}
We see that
\begin{equation*}
\dist (\ell_+,\ell_-)=\frac{2\mu \gamma}{\vert q^\prime -q\vert}
\cdot \cos \theta =\frac{2\mu \gamma}{\sqrt{(q^\prime
-q)^2+(a^\prime -a)^2}} \leq \frac{4\mu}{\sqrt{(q^\prime
-q)^2+(a^\prime -a)^2}} \, .
\end{equation*}

The set whose area defines $A_{P,P^\prime}(Q,\mu)$ is the
intersection of the strip bounded by $\ell_+$ and $\ell_-$ and the
disc $\D_0$, thus
\begin{equation}\label{2.5}
A_{P,P^\prime} (Q,\mu)\leq 2r_0 \dist (\ell_+,\ell_-) \leq
\frac{8\mu r_0}{\sqrt{(q^\prime -q)^2+(a^\prime -a)^2}} \, .
\end{equation}
We also have
\begin{equation}\label{2.6}
A_{P,P^\prime} (Q,\mu)\neq 0 \quad \mbox{\rm only if}\quad \vert
a^\prime q-aq^\prime \vert \leq 2\mu +\vert a^\prime -a\vert
+\vert q^\prime -q\vert .
\end{equation}

\medskip

\begin{lem}\label{L2.1}
Let $\alpha \in (0,1]$. Let $C$ be a compact set in $\R_+$. Then
for all $\eps >0$ and all $\mu \in C$
\begin{equation*}
\frac{1}{Q^2} \sum\limits_{\substack{P\in \Box_{Q^\alpha} \\
P^\prime \in \Box_Q \\ P\neq P^\prime }}
A_{P,P^\prime}(Q,\mu)=O_{C,\D_0,\eps}(Q^{\alpha-1+\eps}).
\end{equation*}
\end{lem}

{\sl Proof.} The estimate \eqref{2.5} reads as
$A_{P,P^\prime}(Q,\mu) =O_{C,\D_0} \big( \frac{1}{\| PP^\prime
\|})$. Combining it with \eqref{2.6} we see that it suffices to
show that
\begin{equation*}
A_Q :=\sum\limits_{\substack{P\in \Box_{Q^\alpha} ,\, P^\prime \in
\Box_Q \\ \vert a^\prime q-aq^\prime \vert \ll_{C,\eps} \|
PP^\prime \|}} \frac{1}{\| PP^\prime \|} \ll_\eps
Q^{\alpha+1+\eps} .
\end{equation*}
Taking $P^{\prime \prime}=(q^{\prime \prime},a^{\prime
\prime})=(q^\prime -q,a^\prime -a)\in \Box_{2Q}$, we gather
\begin{equation*}
\begin{split}
A_Q & \leq \sum\limits_{\substack{P\in \Box_{Q^\alpha},\, O\neq
P^{\prime \prime}\in \Box_{2Q}\\ \vert a^{\prime
\prime}q-aq^{\prime \prime}\vert \ll_C \| OP^{\prime \prime}\|}}
\frac{1}{\| OP^{\prime \prime}\|} \\ & \leq \sum\limits_{O\neq
P^{\prime \prime}\in \Box_{2Q}} \frac{1}{\| OP^{\prime \prime}\|}
\  \# \Big\{ (q,a)\in [-Q^\alpha,Q^\alpha ]^2 \, :\, \vert
a^{\prime \prime} q-aq^{\prime \prime}\vert \ll_C \| OP^{\prime
\prime}\| \Big\} .
\end{split}
\end{equation*}
The two conditions on $(q,a)$ above yield that $(q,a)$ should
belong to the intersection of a strip of width $\ll_C \frac{\|
OP^{\prime \prime} \|}{\| OP^{\prime \prime}\|}=1$  bounded by the
lines $y=\frac{a^{\prime \prime}}{q^{\prime \prime}} \, x\pm
\alpha_C$ with the square $[-Q^\alpha,Q^\alpha ]^2$. The number of
integer lattice points inside this region is of order $O_C
(Q^\alpha)$, thus
\begin{equation*}
A_Q \ll_C Q^\alpha \sum\limits_{O\neq P^{\prime \prime}\in
\Box_{2Q}} \frac{1}{\| OP^{\prime \prime}\|} =Q^\alpha
\sum\limits_{0<m^2+n^2\leq 4Q^2} \frac{1}{\sqrt{m^2+n^2}} \ .
\end{equation*}
Since $r_2(k)=\{ (m,n)\in \Z^2 \, :\, m^2+n^2 =k\}=O_\eps
(k^\eps)$, this gives
\begin{equation*}
\begin{split}
A_Q & \ll_C Q^\alpha \sum\limits_{k=1}^{4Q^2}
\sum\limits_{m^2+n^2=k} \frac{1}{\sqrt{k}} =Q^\alpha
\sum\limits_{k=1}^{4Q^2} \frac{r_2(k)}{\sqrt{k}} \\
& \ll_\eps Q^\alpha \sum\limits_{k=1}^{4Q^2} k^{\eps-\frac{1}{2}}
\ll Q^\alpha (Q^2)^{\eps+\frac{1}{2}} =Q^{\alpha+1+2\eps} ,
\end{split}
\end{equation*}
as desired. \qed

\medskip

\begin{lem}\label{L2.2}
For every compact set $C\subset \R_+$ and every $\eps >0$, there
exist constants $M_1,M_2>0$ such that
\begin{equation*}
\frac{Q^2}{N} \ G_Q \bigg( \frac{\pi\lambda}{2}-M_1
Q^{-\frac{1}{3}}\bigg) -M_2 Q^{-\frac{1}{3}+\eps} \leq R^{(2)}_Q
(\lambda) \leq \frac{Q^2}{N} \ G_Q \bigg( \frac{\pi
\lambda}{2}+M_1 Q^{-\frac{1}{3}}\bigg)+M_2 Q^{-\frac{1}{3}+\eps} .
\end{equation*}
\end{lem}

{\sl Proof.} The trivial estimate
\begin{equation*}
\| P_{(x,y)}P\|=\| OP\|+O_{\D_0}(1)
\end{equation*}
and \eqref{2.1} yield for all $P,P^\prime \in \Box_Q$, $\lambda
\in C$, $(x,y)\in \D_0$, that
\begin{equation*}
\begin{split}
& \beta_{Q,\lambda} \| P_{(x,y)}P\| \cdot \| P_{(x,y)}P^\prime \|
=\beta_{Q,\lambda} \big( \| OP\| \cdot \| OP^\prime \|+O(Q)\big)
\\ & \qquad =\bigg( \frac{\pi \lambda}{2Q^2}+O_C \Big(
\frac{1}{Q^6}\Big)\bigg) \big( \| OP\| \cdot \| OP^\prime \|
+O_{\D_0} (Q)\big) =\frac{\pi \lambda \| OP\| \cdot \| OP^\prime
\|}{2Q^2}+O_{C,\D_0}\Big( \frac{1}{Q}\Big) \\
& \qquad =\frac{\pi\lambda}{2}\ \gamma_{P,P^\prime}(Q)+O_{C,\D_0}
\Big( \frac{1}{Q}\Big) =\gamma_{P,P^\prime}(Q) \bigg(
\frac{\pi\lambda}{2} +O_{C,\D_0} \Big( \frac{1}{Q}\Big) \bigg) \\
& \qquad =\gamma_{P,P^\prime}(Q) \bigg(
\frac{\pi\lambda}{2}+O_{C,\D_0} \Big(
\frac{1}{Q\gamma_{P,P^\prime}(Q)}\Big) \bigg) .
\end{split}
\end{equation*}

We first analyze the case where $\min\{ \| OP\|,\|OP^\prime \|\}
\geq Q^{\frac{2}{3}}$. In this case $\gamma_{P,P^\prime}(Q)\geq
Q^{-\frac{2}{3}}$, and the relation above and the definitions of $
w_{P,P^\prime}$ and $A_{P,P^\prime}$, yield $M_1>0$ such that
\begin{equation*}
A_{P,P^\prime} \bigg( Q,\frac{\pi\lambda}{2}-M_1 Q^{-\frac{1}{3}}
\bigg) \leq w_{P,P^\prime} (Q,\lambda) \leq A_{P,P^\prime} \bigg(
Q,\frac{\pi\lambda}{2}+M_1 Q^{-\frac{1}{3}}\bigg).
\end{equation*}

When $\min\{ \| OP\|,\|OP^\prime \|\} \leq Q^{\frac{2}{3}}$, we
take $\alpha=\frac{1}{3}$ in Lemma \ref{L2.1}. Since
$\beta_{Q,\lambda}\|P_{(x,y)}P\| \cdot \| P_{(x,y)}P^\prime \|
\ll_C \pi\lambda \gamma_{P,P^\prime}(Q)$ as $Q\rightarrow \infty$,
we get
\begin{equation*}
\sum\limits_{\min\{ \| OP\|,\|OP^\prime \|\} \leq Q^{\frac{2}{3}}}
\hspace{-30pt} w_{P,P^\prime}(Q,\lambda) \quad \ll_C
\hspace{-10pt} \sum\limits_{\min\{ \| OP\|,\|OP^\prime \|\} \leq
Q^{\frac{2}{3}}}\hspace{-30pt} A_{P,P^\prime}(Q,\pi \lambda)\quad
\ll_{C,\eps} Q^{-\frac{1}{3}+\eps}.
\end{equation*}
\qed

\section{A formula for $G_Q(\mu)$}
An immediate consequence of \eqref{2.5} and \eqref{2.6} is that
the contribution to $G_Q(\mu)$ of pairs of points $(P,P^\prime)\in
\tilde{\Box}_Q^2$ with $a^\prime =a$ or with $q^\prime =q$ is
negligible. Indeed, we see from \eqref{2.6} that when $a^\prime
=a\neq 0$, the term $A_{P,P^\prime}(Q,\mu)$ is zero unless $\vert
q^\prime -q\vert \leq 2\mu+\frac{1}{\vert a\vert} \leq 2\mu+1$,
thus the total contribution of such points to $G_Q(\mu)$ is
\begin{equation*}
\ll_C \frac{1}{Q^2} \sum\limits_{\substack{\vert q\vert \leq Q
\\ 0<\vert q^\prime -q\vert \leq 2\mu +1 \\ 0<\vert a\vert \leq Q}}
\frac{8\mu r_0}{\vert q^\prime -q\vert} \ll_C \frac{1}{Q}
\sum\limits_{\substack{\vert q\vert \leq Q \\ 0<\vert q^\prime
-q\vert<2\mu +1}} \frac{1}{\vert q^\prime -q\vert} \ll_C \frac{\ln
Q}{Q} \ .
\end{equation*}
The contribution of pairs of points $(P,P^\prime)\in
\tilde{\Box}_Q^2$ with $a^\prime =a=0$ to $G_Q (\mu)$ is
\begin{equation*}
\ll_C \frac{1}{Q^2} \sum\limits_{\substack{\vert q\vert,\vert
q^\prime \vert \leq Q \\ q^\prime \neq q}} \frac{1}{\vert q^\prime
-q\vert} \ll \frac{\ln Q}{Q}\ .
\end{equation*}
Similar estimates in the case $q^\prime =q$ show that
\begin{equation}\label{2.9}
G_Q(\mu)=\frac{1}{Q^2} \sum\limits_{\substack{(P,P^\prime) \in
\tilde{\Box}_Q^2
\\ a^\prime \neq a,\,q^\prime \neq q}} A_{P,P^\prime}(Q,\mu)\
+O_{C,\D_0} \bigg( \frac{\ln Q}{Q} \bigg) .
\end{equation}
As a result, we shall subsequently assume that $a^\prime \neq a$
and $q^\prime \neq q$. We now set
\begin{equation*}
\alpha =\frac{a^\prime -a}{q^\prime -q} \, ;\qquad \beta
=\frac{aq^\prime -a^\prime q}{q^\prime -q} \, ; \qquad \gamma_0
=\frac{\gamma}{\vert q^\prime -q\vert}=\frac{\sqrt{q^2+a^2}
\sqrt{q^{\prime 2}+a^{\prime 2}}}{Q^2 \vert q^\prime -q\vert} \, .
\end{equation*}

The remainder of this section is elementary and is concerned with
putting $G_Q(\mu)$ in a tidy form, suitable for a precise
estimation which will be completed in the next section.

Let $C_0$ denote the center of $\D_0$, let $\ell^\prime$ be the
line passing through $C_0$ and perpendicular to $\ell$, and denote
by $A_+$ and $A_-$ the intersections of $\ell^\prime$ with the
circle $\partial \D_0$, by $E_0$ the intersection of $\ell^\prime$
and $\ell$, and by $E_\pm$ the intersection of $\ell^\prime$ with
$\ell_\pm$. Direct computation gives
\begin{equation*}
\begin{split}
& x_{A_\pm} =x_0 \mp \frac{\alpha r_0}{\sqrt{\alpha^2+1}} =x_0 \mp
\frac{(a^\prime -a)r_0}{\sqrt{(q^\prime -q)^2+(a^\prime -a)^2}}\,
; \qquad x_{E_\pm} =\frac{\alpha y_0+x_0-\alpha \beta\mp \alpha
\mu \gamma_0}{\alpha^2+1}\ ; \\
& \| E_- E_+\| =\dist(\ell_+,\ell_-)=\frac{\vert x_{E_-} -x_{E_+}
\vert}{\vert \sin \theta \vert}=\frac{2\mu
\gamma_0}{\sqrt{\alpha^2+1}} = \frac{\mu \sqrt{q^2+a^2}
\sqrt{q^{\prime 2}+a^{\prime 2}}}{Q^2 \| PP^\prime\|} \ .
\end{split}
\end{equation*}
\begin{figure}[ht]
\includegraphics*[scale=0.7, bb=0 0 320 340]{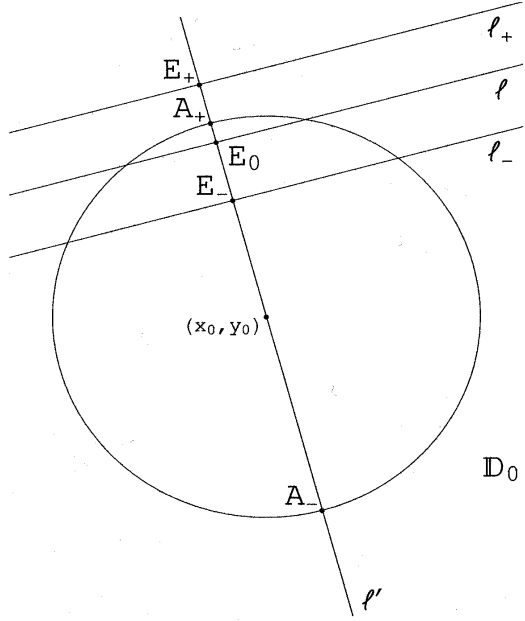}
\caption{The intersection of the strip bounded by $\ell_+$ and
$\ell_-$ with the disc $\D_0$.}
\end{figure}

While ordering the points $x_{E_+}<x_{E_-}$ and $x_{A_+}<x_{A_-}$
the following situations may occur:

{\bf Case 1.} \quad $x_{E_+} <x_{A_+}<x_{A_-} <x_{E_-}$, that is
\begin{equation*}
\frac{\alpha y_0+x_0-\alpha \beta-\alpha \mu \gamma_0}{\alpha^2+1}
<x_0-\frac{\alpha r_0}{\sqrt{\alpha^2+1}} < x_0+\frac{\alpha
r_0}{\sqrt{\alpha^2+1}} <\frac{\alpha y_0+x_0-\alpha \beta +\alpha
\mu \gamma_0}{\alpha^2+1} \ .
\end{equation*}
This gives $\mu \gamma_0 >r_0 \sqrt{\alpha^2+1}$, hence
\begin{equation*}
r_0 \sqrt{(q^\prime -q)^2+(a^\prime -a)^2} <\frac{\mu
\sqrt{q^2+a^2} \sqrt{q^{\prime 2}+a^{\prime 2}}}{Q^2} \leq 2\mu .
\end{equation*}
Suppose first that $\vert a^\prime -a\vert \leq \vert q^\prime
-q\vert$. By \eqref{2.5} we know that for fixed $(q,q^\prime)$,
the expression $D=aq^\prime-a^\prime q$ will only take values
between $-2\mu -\frac{4\mu}{r_0}$ and $2\mu +\frac{4\mu}{r_0}$.
Hence the number of solutions $(a,a^\prime)$ of
$aq^\prime-a^\prime q=D$ is of order $O_C (d)$, where $d$ is the
greatest common divisor of $q$ and $q^\prime$. But $d\leq
\frac{2\mu}{r_0}$, hence this is actually $O_C (1)$ and the
contribution to $G_Q$ is
\begin{equation*}
\ll_C \frac{1}{Q^2} \sum\limits_{\substack{\vert q\vert,\vert
q^\prime \vert \leq Q \\ 0<\vert q^\prime -q \vert \ll 1}}
\sum\limits_{\substack{\vert a\vert , \vert a^\prime \vert \leq Q \\
A_{P,P^\prime} \neq 0 \\ 0<\vert a^\prime -a\vert \leq \vert
q^\prime -q\vert}} 1\ \ll \ \frac{1}{Q^2}
\sum\limits_{\substack{\vert q\vert ,\vert q^\prime \vert \leq Q \\
0< \vert q^\prime -q\vert \ll 1}} 1 \ \ll \ \frac{1}{Q} \ .
\end{equation*}
The case $\vert q^\prime -q\vert \leq \vert a^\prime -a\vert$ is
settled similarly by first summing over $(a,a^\prime)$.

{\bf Case 2.} \quad $x_{A_+}<x_{E_+}<x_{A_-}<x_{E_-}$, that is
\begin{equation*}
\vert \mu \gamma_0-r_0 \sqrt{\alpha^2+1} \vert
=r_0\sqrt{\alpha^2+1}-\mu \gamma_0 <y_0 -\alpha x_0 -\beta <\mu
\gamma_0 +r_0 \sqrt{\alpha^2+1}\ ,
\end{equation*}
or equivalently
\begin{equation*}
\big| a^\prime q-aq^\prime +(q^\prime -q)y_0-(a^\prime -a)x_0 -r_0
 \sqrt{(q^\prime -q)^2+(a^\prime -a)^2} \, \big| <\mu \gamma .
\end{equation*}
The change of variables $a^\prime -a=a^{\prime \prime}$, $q^\prime
-q=q^{\prime \prime}$ gives
\begin{equation*}
\vert a^{\prime \prime} q-aq^{\prime \prime} -r_0 \sqrt{q^{\prime
\prime 2}+a^{\prime \prime 2}} +q^{\prime \prime} y_0 -a^{\prime
\prime} x_0 \vert <\mu \gamma \leq 4\mu .
\end{equation*}
So, keeping $a^{\prime \prime}$ and $q^{\prime \prime}$ fixed, the
range of $a^{\prime \prime} q-aq^{\prime \prime}$ has cardinality
$O_C (1)$. Now the equation $a^{\prime \prime} q-aq^{\prime
\prime}=K$ has either no solution $(q,a)$ when $d=\gcd (a^{\prime
\prime},q^{\prime \prime})$ does not divide $K$, or has $O\big(
\frac{dQ}{q^{\prime \prime}}\big)$ solutions $(q,a)$ when $d$
divides $K$. Thus the contribution of terms $A_{P,P^\prime}$ with
$q^{\prime \prime 2}+a^{\prime \prime 2}=(q^\prime -q)^2+(a^\prime
-a)^2 >Q$ is
\begin{equation*}
\ll_C \frac{1}{Q^2} \sum\limits_{d=1}^Q
\sum\limits_{\substack{0<\vert q_0^{\prime \prime}\vert ,\vert
a_0^{\prime \prime} \vert \leq \big[ \frac{Q}{d}\big]
\\ \gcd (q_0^{\prime \prime},a_0^{\prime \prime})=1}}
\frac{Q}{q_0^{\prime \prime}} \cdot \frac{1}{\sqrt{Q}} \
 \leq \  \frac{1}{Q\sqrt{Q}} \sum\limits_{d=1}^Q \frac{Q}{d}
\sum\limits_{0<\vert q_0^{\prime \prime}\vert \leq \big[
\frac{Q}{d}\big]} \frac{1}{q_0^{\prime \prime}} \  \ll \
\frac{\ln^2 Q}{\sqrt{Q}} \ \ll_\delta Q^{-\frac{1}{2}+\delta} .
\end{equation*}
The contribution of terms $A_{P,P^\prime}$ with $q^{\prime \prime
2}+a^{\prime \prime 2}=(q^\prime -q)^2+(a^\prime -a)^2 \leq Q$ is
\begin{equation*}
\ll_C \frac{1}{Q^2} \sum\limits_{1\leq d\leq \sqrt{Q}}
\sum\limits_{0<\vert q_0^{\prime \prime}\vert,\vert a_0^{\prime
\prime}\vert \ll \big[ \frac{\sqrt{Q}}{d}\big]}
\frac{Q}{q_0^{\prime \prime}} \cdot \frac{1}{dq_0^{\prime \prime}}
\ \leq \ \frac{1}{Q} \sum\limits_{d,q_0^{\prime \prime}=1}^\infty
\frac{\sqrt{Q}}{d^2 q_0^{\prime \prime}} \ \ll \ Q^{-\frac{1}{2}}
\ .
\end{equation*}

{\bf Case 3.} \quad $x_{E_+}<x_{A_+}<x_{E_-}<x_{A_-}$, that is
\begin{equation*}
-\mu \gamma_0 -r_0 \sqrt{\alpha^2+1} <y_0 -\alpha x_0 -\beta
<-\vert \mu \gamma_0 -r_0 \sqrt{\alpha^2+1} \vert .
\end{equation*}

We infer as in Case 2 that the contribution of $A_{P,P^\prime}$ is
in this case too $O_\delta (Q^{-\frac{1}{2}+\delta})$.

{\bf Case 4.} \quad $x_{A_+}<x_{E_+}<x_{E_-}<x_{A_-}$, that is
\begin{equation*}
\mu \gamma_0 -r_0 \sqrt{\alpha^2+1} <y_0-\alpha x_0-\beta <-\mu
\gamma_0 +r_0 \sqrt{\alpha^2+1}\ ,
\end{equation*}
or equivalently
\begin{equation*}
\vert L_{P,P^\prime}(x_0,y_0)\vert <r_0 \sqrt{(q^\prime
-q)^2+(a^\prime -a)^2}-\mu\gamma=r_0 \| PP^\prime \|-\mu \gamma.
\end{equation*}
Denote $k=q^\prime -q$ and $\ell=a^\prime -a$. The interval
$I_{k,\ell}=[r_0 \sqrt{k^2+\ell^2} -ky_0+\ell x_0 -\mu \gamma, r_0
\sqrt{k^2+\ell^2} -ky_0+\ell x_0 ]$ has length $\mu \gamma \ll_C
1$. Hence we find that the contribution of terms $A_{P,P^\prime}$
for which $r_0 \|PP^\prime \| -\mu \gamma <L_{P,P^\prime}
(x_0,y_0)<r_0 \| PP^\prime \|$ to $G_Q$ is
\begin{equation*}
\begin{split}
& \ll_C \frac{1}{Q^2} \sum\limits_{0<\vert k\vert,\vert \ell \vert
\leq Q} \sum\limits_{\substack{\vert q^\prime \vert ,\vert
a^\prime \vert \leq Q  \\ -ka^\prime -\ell q^\prime \in
I_{k,\ell}}} \frac{1}{\sqrt{k^2+\ell^2}} \ll \frac{1}{Q^2}
\sum\limits_{k,\ell=1}^Q \frac{Q\gcd (k,l)}{k}\cdot
\frac{1}{\sqrt{k^2+\ell^2}}
\\ & \leq \frac{1}{Q^2} \sum\limits_{d=1}^Q \frac{1}{d}
\sum\limits_{k_0,\ell_0=1}^{\big[ \frac{Q}{d}\big]} \frac{1}{k_0
\sqrt{k_0^2+\ell_0^2}} \ll \frac{1}{Q} \sum\limits_{d=1}^Q
\frac{1}{d} \sum\limits_{k_0=1}^{\big[ \frac{Q}{d}\big]}
\frac{1}{k_0} \sum\limits_{\ell_0=1}^{\big[
\frac{Q}{d}\big]} \frac{1}{\ell_0} \\
& \ll \frac{\ln^3 Q}{Q} \ .
\end{split}
\end{equation*}
One shows similarly that the contribution of points $P,P^\prime$
for which $-r_0\| PP^\prime \|+\mu \gamma
<L_{P,P^\prime}(x_0,y_0)<r_0 \| PP^\prime\|$ is of the same order.
Therefore by \eqref{2.9} and the previous considerations we infer
that
\begin{equation}\label{2.10}
G_Q (\mu) =\frac{1}{Q^2} \hspace{-20pt}
\sum\limits_{\substack{(P,P^\prime) \in \tilde{\Box}_Q ^2\\
a^\prime \neq a,\, q^\prime \neq q \\ \vert
L_{P,P^\prime}(x_0,y_0)\vert <r_0 \| PP^\prime \|}} \hspace{-25pt}
A_{P,P^\prime} (Q,\mu) +O_{C,\D_0,\delta}
(Q^{-\frac{1}{2}+\delta}).
\end{equation}

\begin{figure}[ht]
\includegraphics*[scale=0.7, bb=0 0 250 250]{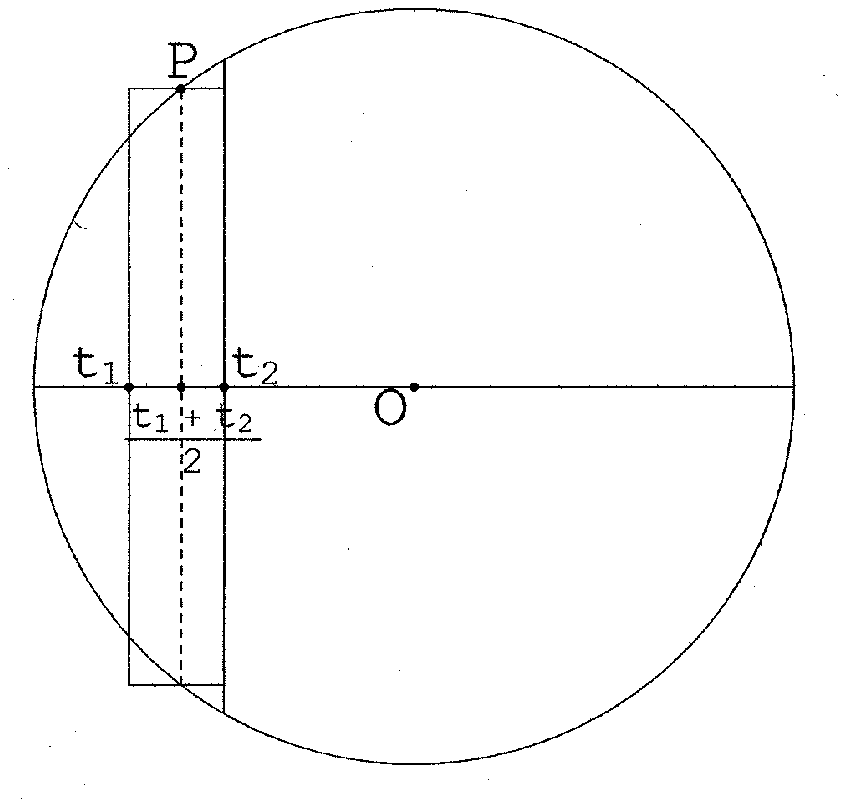}
\caption{}
\end{figure}

Next $S_{P,P^\prime}$ is approximated by elementary calculus.

\medskip

\begin{lem}\label{LII1}
The area of the region inside the circle of radius $r_0$ centered
at the origin and inside the strip bounded by the vertical lines
$y=t_1$ and $y=t_2$, $-r_0<t_1\leq t_2 <r_0$, is given, for small
$t_2-t_1$, by
\begin{equation*}
A_{r_0}(t_1,t_2)=2(t_2-t_1) \sqrt{r_0^2-\bigg(
\frac{t_1+t_2}{2}\bigg)^2 }\, + O_{r_0}\big(
(t_2-t_1)^{\frac{3}{2}}\big)
\end{equation*}
\end{lem}

{\sl Proof.} The error is seen to be given (see Figure 2) by
\begin{equation*}
\int\limits_{t_1}^{\frac{t_1+t_2}{2}} \Bigg( \sqrt{r_0^2 -\bigg(
\frac{t_1+t_2}{2} \bigg)^2} \ -\sqrt{r_0^2-t^2} \ \Bigg)
dt+\int\limits_{\frac{t_1+t_2}{2}}^{t_2} \Bigg( \sqrt{r_0^2-t^2}
-\sqrt{r_0^2 -\bigg( \frac{t_1+t_2}{2} \bigg)^2} \ \Bigg) dt.
\end{equation*}
It is $\ll (t_2-t_1)^{3/2}$ as a result of
\begin{equation*}
\big| \sqrt{r_0^2-x^2} -\sqrt{r_0^2 -y^2} \, \big| \leq
\sqrt{\vert x^2-y^2 \vert} \leq 2\sqrt{\vert x-y\vert} .
\end{equation*}
\qed

We take
\begin{equation*}
t_{E_\pm} :=\frac{x_{E_\pm}-x_{E_0}}{\sin \theta}
=\frac{\sqrt{1+\alpha^2}}{\alpha} \  (x_{E_\pm}-x_{E_0})
=\frac{L_{P,P^\prime}(x_0,y_0)\mp \mu \gamma}{\| PP^\prime \|} \ .
\end{equation*}

Notice now that
\begin{equation}\label{2.11}
t_{E_-}-t_{E_+} =\frac{2\mu \gamma}{\| PP^\prime \|} \ll_C
\frac{1}{\| PP^\prime \|} \ ,
\end{equation}
denote $k=q^\prime -q$, $\ell=a^\prime -a$,
\begin{equation}\label{2.12}
J_{k,\ell}=[-ky_0+\ell x_0-r_0 \sqrt{k^2+\ell^2},-ky_0+\ell
x_0+r_0 \sqrt{k^2+\ell^2}] ,
\end{equation}
and find that the contribution of the error provided by Lemma
\ref{LII1} in \eqref{2.10} is
\begin{equation}\label{2.13}
\begin{split} \frac{1}{Q^2} \hspace{-20pt}
\sum\limits_{\substack{(P,P^\prime) \in \tilde{\Box}_Q^2 \\
q^\prime \neq q,\, a^\prime \neq a \\ \vert
L_{P,P^\prime}(x_0,y_0)\vert <r_0 \| PP^\prime \|}} \hspace{-20pt}
\frac{1}{\| PP^\prime \|^{3/2}} & =\frac{1}{Q^2}
\sum\limits_{0<\vert k\vert,\vert \ell \vert \leq Q}
\sum\limits_{\substack{\vert q^\prime \vert ,\vert a^\prime \vert
\leq Q \\ ka^\prime -\ell q^\prime \in J_{k,\ell}}}
\frac{1}{(k^2+\ell^2)^{3/4}}
\\ & \ll \frac{1}{Q^2} \sum\limits_{k,\ell=1}^Q \frac{Q\gcd
(k,\ell)}{k} \cdot \sqrt{k^2+\ell^2} \cdot
\frac{1}{(k^2+\ell^2)^{3/4}} \\ & \leq \frac{1}{Q^2}
\sum\limits_{d=1}^Q \sum\limits_{k_0,\ell_0=1}^{\big[
\frac{Q}{d}\big]} \frac{Qd}{dk_0} \cdot
d(k_0^2+\ell_0^2)^{\frac{1}{2}} \cdot \frac{1}{d^{3/2}
(k_0^2+\ell_0^2)^{3/4}} \\ & =\frac{1}{Q} \sum\limits_{d=1}^Q
\frac{1}{d^{1/2}} \sum\limits_{k_0,\ell_0=1}^{\big[
\frac{Q}{d}\big]} \frac{1}{k_0(k_0^2+\ell_0^2)^{1/4}} \\
& \leq \frac{1}{Q} \sum\limits_{d=1}^Q \frac{1}{d^{1/2}}
\sum\limits_{k_0}^{\big[ \frac{Q}{d}\big]} \frac{1}{k_0}
\sum\limits_{\ell_0=1}^{\big[ \frac{Q}{d}\big]}
\frac{1}{\ell_0^{1/2}} \\
& \ll \frac{\ln Q}{Q} \sum\limits_{d=1}^Q \frac{1}{d^{1/2}} \cdot
\frac{Q^{1/2}}{d^{1/2}} \\
& \ll \frac{\ln^2 Q}{Q^{1/2}} \ .
\end{split}
\end{equation}
Therefore by \eqref{2.10}, Lemma \ref{LII1}, \eqref{2.11} and
\eqref{2.13} we find that
\begin{equation}\label{2.14}
G_Q(\mu) =\frac{1}{Q^2} \hspace{-20pt}
\sum\limits_{\substack{(P,P^\prime) \in \tilde{\Box}_Q^2 \\
q^\prime \neq q,\, a^\prime \neq a \\ \vert
L_{P,P^\prime}(x_0,y_0)\vert <r_0 \| PP^\prime \|}} \hspace{-25pt}
B_{P,P^\prime}(\mu) \ +O_{C,\D_0,\delta}
(Q^{-\frac{1}{2}+\delta}),
\end{equation}
where $B_{P,P^\prime}(Q,\mu)$ denotes the contribution of the main
term in Lemma \ref{LII1} to \eqref{2.10}, that is
\begin{equation}\label{2.15}
\begin{split}
B_{P,P^\prime}(Q,\mu) & =2(t_{E_-}-t_{E_+})\sqrt{r_0^2-\bigg(
\frac{t_{E_-}+t_{E_+}}{2} \bigg)^2}  =2\cdot \frac{2\mu \gamma}{\|
PP^\prime \|} \sqrt{r_0^2-\frac{L_{P,P^\prime}(x_0,y_0)^2}{\|
PP^\prime \|^2}}
\\ & =\frac{4\mu \gamma\sqrt{r_0^2
\| PP^\prime \|^2 -L_{P,P^\prime} (x_0,y_0)^2}}{\| PP^\prime\|^2}
\ .
\end{split}
\end{equation}

Finally we show that, if $a^\prime -a\leq q^\prime -q$, one can
replace $\gamma$ by $\frac{1}{Q^2} \,qq^\prime \big(
1+\frac{(a^\prime -a)^2}{(q^\prime -q)^2}\big)$ in \eqref{2.15}
and \eqref{2.14}. Since $\vert L_{P,P^\prime}(x_0,y_0)\vert <r_0
\| PP^\prime \|$, then $\vert a^\prime q-aq^\prime \vert \leq
\big(r_0+\sqrt{x_0^2+y_0^2}\, \big) \| PP^\prime\|$, and
\begin{equation*}
\Bigg| q\sqrt{1+\frac{(a^\prime -a)^2}{(q^\prime -q)^2}}\
-q\sqrt{1+\frac{a^2}{q^2}}\ \bigg| \ll q\Bigg| \frac{a^\prime
-a}{q^\prime -q} -\frac{a}{q} \bigg| =\frac{\vert a^\prime
q-aq^\prime \vert}{q^\prime -q} \ll \frac{\vert a^\prime
q-aq^\prime \vert}{\| PP^\prime \|} \ll_{\D_0} 1.
\end{equation*}
This gives
\begin{equation*}
\sqrt{q^2+a^2}=q\sqrt{1+\frac{(a^\prime -a)^2}{(q^\prime -q)^2}} \
+O_{\D_0}(1),
\end{equation*}
and similarly
\begin{equation*}
\sqrt{q^{\prime 2}+a^{\prime 2}}=q^\prime \sqrt{1+\frac{(a^\prime
-a)^2}{(q^\prime -q)^2}} \ +O_{\D_0}(1).
\end{equation*}
Hence one can replace $B_{P,P^\prime} (Q,\mu)$ in \eqref{2.14} by
\begin{equation*}
W_{P,P^\prime} (Q,\mu) =\frac{4\mu qq^\prime \sqrt{r_0^2 \|
PP^\prime \|^2 -L_{P,P^\prime}(x_0,y_0)^2}}{Q^2 \max\{
(q^\prime-q)^2,(a^\prime -a)^2\}}  \ ,
\end{equation*}
at the cost of an error which is
\begin{equation*}
\begin{split}
\ll & \frac{1}{Q^2} \hspace{-20pt}
\sum\limits_{\substack{(P,P^\prime) \in \tilde{\Box}_Q^2
\\ q^\prime \neq q,\ a^\prime \neq a \\
\vert L_{P,P^\prime}(x_0,y_0)\vert <r_0 \| PP^\prime \|}}
\hspace{-25pt} \frac{1}{Q\sqrt{(q^\prime -q)^2+(a^\prime -a)^2}} \\
& \ll \frac{1}{Q^3} \sum\limits_{0<\vert k\vert,\vert \ell \vert
\leq Q} \sqrt{k^2+\ell^2} \cdot \frac{Q\gcd (k,\ell)}{k} \cdot
\frac{1}{\sqrt{k^2+\ell^2}} \\
& \leq \frac{1}{Q^2} \sum\limits_{d=1}^Q \
\sum\limits_{k_0,\ell_0=1}^{\big[ \frac{Q}{d}\big]} \frac{Q}{k_0}
\ll \frac{1}{Q^2} \sum\limits_{d=1}^Q \frac{Q}{d}\cdot \ln Q
=\frac{\ln^2 Q}{Q} \ .
\end{split}
\end{equation*}

In summary, we have shown for any $\mu$ in a fixed compact set
$C\subset \R_+$, that
\begin{equation}\label{2.16}
G_Q(\mu)=\frac{4\mu}{Q^4} \hspace{-30pt}
\sum\limits_{\substack{(P,P^\prime) \in \tilde{\Box}_Q^2 \\
q^\prime \neq q,\ a^\prime \neq a \\ r_0 \| PP^\prime \| \geq
\vert L_{P,P^\prime}(x_0,y_0)\vert}} \hspace{-32pt}
\frac{qq^\prime \sqrt{r_0^2 \| PP^\prime
\|^2-L_{P,P^\prime}(x_0,y_0)^2}}{\max \{ (q^\prime -q)^2,(a^\prime
-a)^2\}} +O_{C,\D_0,\delta} (Q^{-\frac{1}{2}+\delta})
\end{equation}
as $Q\rightarrow \infty$.

\bigskip

\section{Estimating the sum $S_Q$}

By reflecting $\D_0$ about the axes and about the line $y=x$, we
see that it suffices to only estimate the contribution $A_Q (\mu)$
to $G_Q (\mu)$ of points $(P,P^\prime)\in (0,Q]^2$ with $0<\alpha
=\frac{a^\prime -a}{q^\prime -q}\leq 1$. We thus consider
\begin{equation*}
A_Q(\mu) =4\mu S_Q ,
\end{equation*}
where
\begin{equation*}
\begin{split}
S_Q & =\frac{1}{Q^4} \hspace{-20pt}
\sum\limits_{\substack{0<q,q^\prime \leq Q\\
0<a,a^\prime \leq Q \\ r_0 \| PP^\prime \| >\vert L_{P,P^\prime}(x_0,y_0)\vert \\
0<\frac{a^\prime -a}{q^\prime -q}\leq 1}} \hspace{-20pt}
\frac{qq^\prime \sqrt{r_0^2 \| PP^\prime
\|^2-L_{P,P^\prime}(x_0,y_0)^2}}{(q^\prime -q)^2} \\ &
=\frac{2}{Q^4} \hspace{-20pt}\sum\limits_{\substack{0<a<a^\prime \leq Q \\
0<q<q^\prime \leq Q \\ r_0 \| PP^\prime \|
>\vert L_{P,P^\prime}(x_0,y_0)\vert \\ 0<\frac{a^\prime
-a}{q^\prime -q} \leq 1}} \hspace{-20pt}\frac{qq^\prime
\sqrt{r_0^2 \big( (q^\prime -q)^2+(a^\prime
-a)^2\big)-L_{P,P^\prime}(x_0,y_0)^2}}{(q^\prime -q)^2} \ .
\end{split}
\end{equation*}
Then we gather from \eqref{2.16} and the above formula for $S_Q$
that
\begin{equation}\label{3.1}
G_Q(\mu)=8A_Q(\mu) +O_{C,\D_0,\delta}(
Q^{-\frac{1}{10}+\delta})=32\mu S_Q+O_{C,\D_0,\delta}(
Q^{-\frac{1}{10}+\delta}).
\end{equation}
Changing $q$ to $q^\prime -q$ and $a$ to $a^\prime -a$, we may
write
\begin{equation*}
S_Q=\frac{2}{Q^4} \hspace{-36pt}
\sum\limits_{\substack{0<a<a^\prime \leq Q \\
0<q<q^\prime \leq Q \\ a\leq q \\
r_0\sqrt{q^2+a^2} >\vert y_0 q-x_0 a+aq^\prime -a^\prime q\vert}}
\hspace{-40pt}  \frac{(q^\prime -q)q^\prime
\sqrt{r_0^2(q^2+a^2)-(y_0 q-x_0 a+aq^\prime -a^\prime q)^2}}{q^2}\
.
\end{equation*}
Putting
\begin{equation*}
D=aq^\prime -a^\prime q
\end{equation*}
and taking $J_{q,a}$ as in \eqref{2.12}, that is
\begin{equation*}
J_{q,a}=[-qy_0 +ax_0 -r_0 \sqrt{q^2+a^2},-qy_0 +ax_0 +r_0
\sqrt{q^2+a^2} \, ],
\end{equation*}
we get
\begin{equation}\label{3.2}
S_Q =\frac{2}{Q^4}
\sum\limits_{1\leq a\leq q\leq Q} \ \sum\limits_{D\in J_{q,a}}
\sum\limits_{\substack{q^\prime \in [q, Q]
\\ a^\prime \in [a, Q] \\ aq^\prime -a^\prime q=D}}
\frac{(q^\prime -q)q^\prime \sqrt{r_0^2(q^2+a^2)-(y_0 q-x_0
a+D)^2}}{q^2} \ .
\end{equation}

We will prove the following result.

\medskip

\begin{prop}\label{P3.1}
$ \quad \displaystyle S_Q=\frac{\pi r_0^2}{6}+O_{\D_0,\delta}
(Q^{-\frac{1}{10}+\delta})\qquad \mbox{for all $\delta
>0$.}$\end{prop}

From this and \eqref{3.1} we infer the following

\medskip

\begin{cor}\label{C3.2}
$\quad \displaystyle G_Q(\mu)=\frac{16\pi r_0^2
\mu}{3}+O_{C,\D_0,\delta}(Q^{-\frac{1}{10}+\delta}).$
\end{cor}

Theorem \ref{T1.1} now follows combining Corollary \ref{C3.2} with
Lemma \ref{L2.2}.

\bigskip

We now start the proof of Proposition \ref{P3.1}. We first lay out
some notation and prove an elementary calculus lemma. Fix
$\alpha_0,\beta_0 \in \R$ and consider the function
\begin{equation*}
\Phi(t,x)=\Phi_{\alpha_0,\beta_0}(t,x)=1+t^2
-(\beta_0-t\alpha_0+x)^2,
\end{equation*}
and the domain
\begin{equation*}
\DD=\DD_{\alpha_0 ,\beta_0}=\{ (t,x)\, :\, 0\leq t\leq 1,\
\Phi(t,x)\geq 0\}.
\end{equation*}
\begin{figure}[ht]
\includegraphics*[scale=0.65, bb=0 0 320 280]{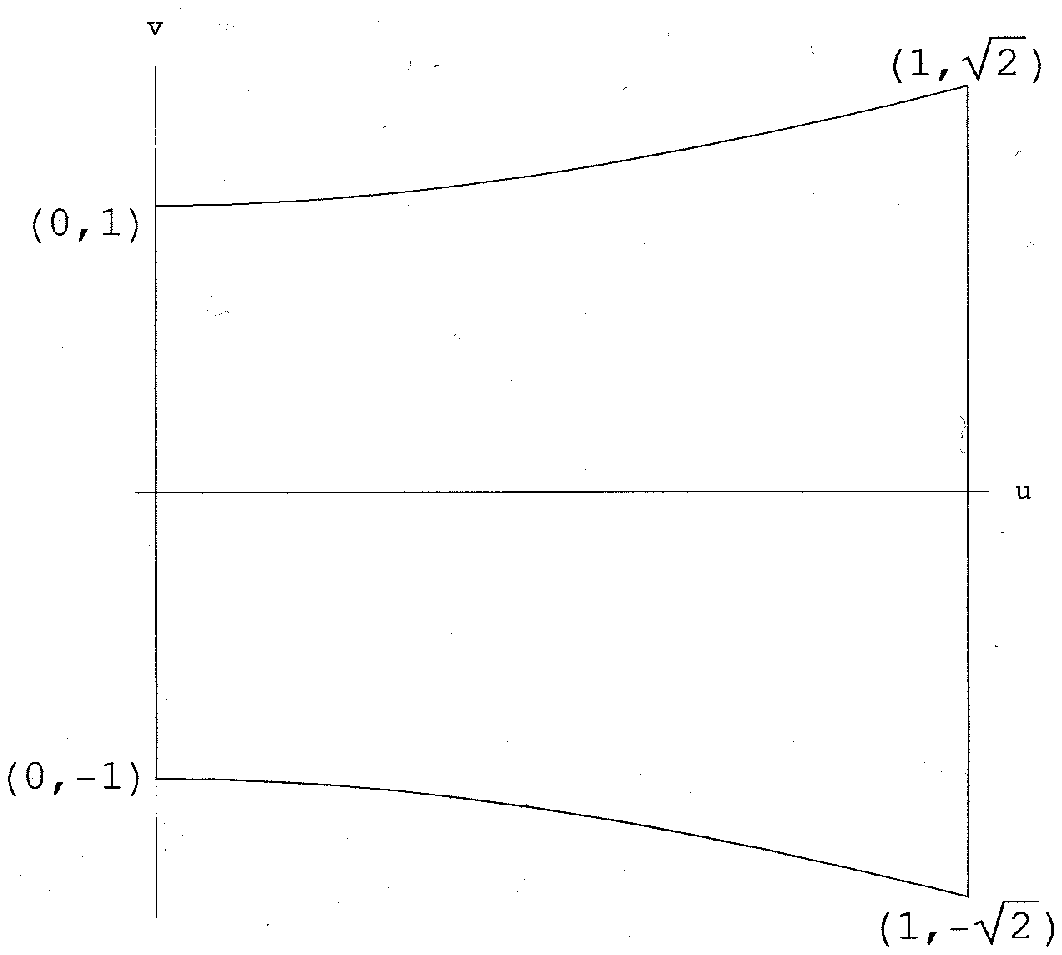}
\includegraphics*[scale=0.65, bb=0 0 300 280]{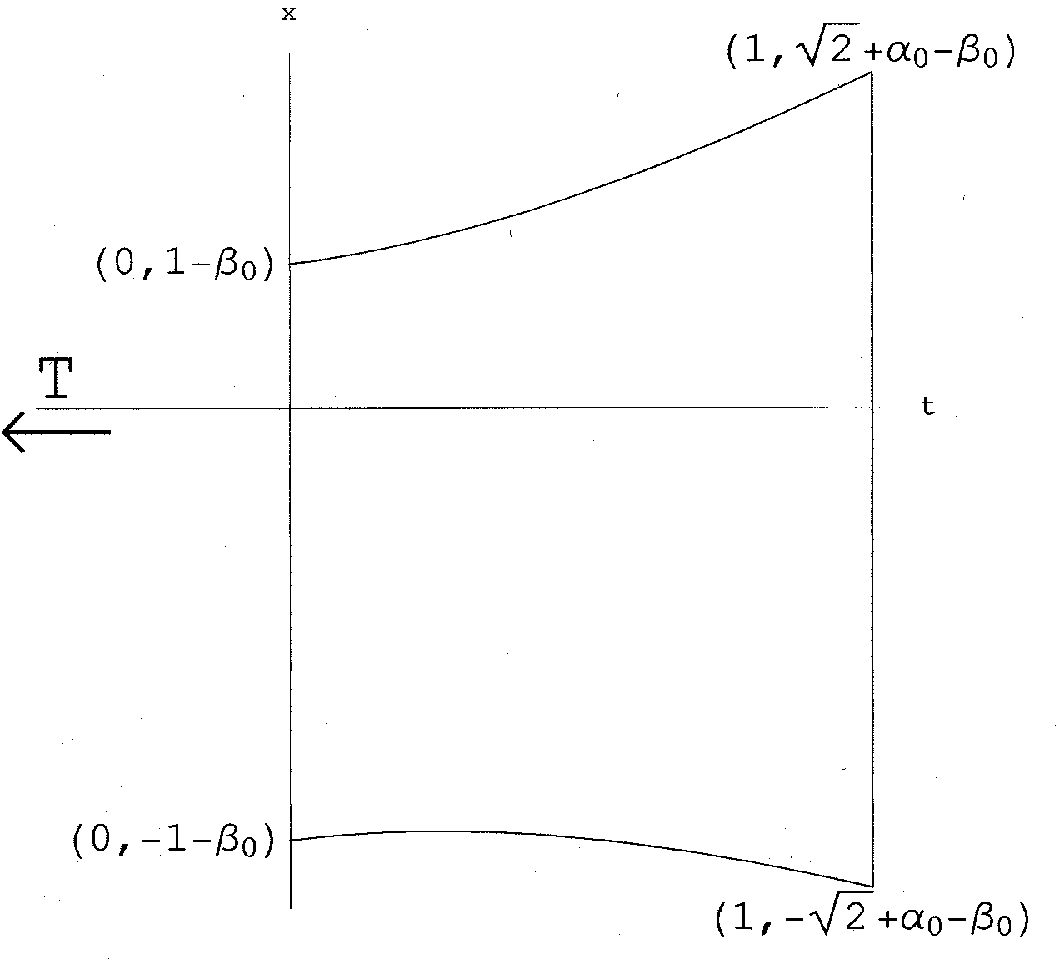}
\caption{The regions $\DD_{0,0}$ and $\DD_{\alpha_0,\beta_0}$.}
\end{figure}

Consider also the projection $\mathrm{pr}_2 \DD$ of $\DD$ on the
second coordinate, the $x$-section
\begin{equation}\label{3.3}
I_x =\{ t\in [0,1]\, :\, \Phi(t,x)\geq 0\},
\end{equation}
and the $t$-section
\begin{equation}\label{3.4}
J_t =\{ x\in \mathrm{pr}_2 \DD \, :\, \Phi(t,x)\geq 0\}.
\end{equation}
Define the function $\psi =\psi_{\alpha_0,\beta_0}:\mathrm{pr}_2
\DD \rightarrow [0,\infty)$ by
\begin{equation*}
\psi (x)=\int\limits_{I_x} \sqrt{\Phi (t,x)}\ dt.
\end{equation*}

\medskip

\begin{lem}\label{L3.0}
For every $\alpha_0,\beta_0 \in \R$ one has

\medskip

\mbox{\rm (i)} \quad $\displaystyle \int\limits_{\mathrm{pr}_2
\DD_{\alpha_0,\beta_0}} \hspace{-15pt} \vert \psi^\prime
(x)\vert\, dx \leq \sqrt{2}+\ln (1+\sqrt{2}).$

\bigskip

\mbox{\rm (ii)} \quad $\displaystyle \int\limits_{\mathrm{pr}_2
\DD_{\alpha_0,\beta_0}} \hspace{-10pt} \psi(x)\, dx=
\iint\limits_{\DD_{\alpha_0,\beta_0}} \sqrt{\Phi (t,x)}\ dt\,
dx=\frac{2\pi}{3} \ .$
\end{lem}

{\sl Proof.} (i) By the definition of $\Phi$ it is seen that $I_x$
is the union of one or two intervals $[a(x),b(x)]$, where $a(x)$
and $b(x)$ are equal to $0$, $1$, or a root of $\Phi(t,x)=0$. In
all these cases
\begin{equation*}
\Phi \big( a(x),x\big) a^\prime (x)=\Phi \big( b(x),x\big)b^\prime
(x)=0,
\end{equation*}
and as a result we get
\begin{equation*}
\frac{d}{dx} \int\limits_{a(x)}^{b(x)} \sqrt{\Phi(t,x)}\ dt=
\int\limits_{a(x)}^{b(x)} \frac{t\alpha_0-\beta_0-x}{\sqrt{\Phi
(t,x)}}\ dt
\end{equation*}
and
\begin{equation*}
\psi^\prime (x)=\frac{d}{dx} \int\limits_{I_x} \sqrt{\Phi(t,x)}\
dt= \int\limits_{I_x } \frac{t\alpha_0-\beta_0-x}{\sqrt{\Phi
(t,x)}}\ dt.
\end{equation*}
Using the triangle inequality and the change of variables
$(u,v)=T(t,x)=(t,t\alpha_0-\beta_0 -x)$ we obtain
\begin{equation*}
\begin{split}
\int\limits_{\mathrm{pr}_2 \DD_{\alpha_0,\beta_0}} \hspace{-10pt}
\vert \psi^\prime (x)\vert\, dx & =\int\limits_{\mathrm{pr}_2
\DD_{\alpha_0,\beta_0}} \left|\ \int\limits_{I_x }
\frac{t\alpha_0-\beta_0-x}{\sqrt{\Phi (t,x)}}\ dt \right|\, dx
\leq \int\limits_{\mathrm{pr}_2 \DD_{\alpha_0,\beta_0}}
\int\limits_{I_x} \frac{\vert
t\alpha_0-\beta_0-x\vert}{\sqrt{\Phi(t,x)}}\ dt \, dx
\\ & =\iint\limits_{\DD_{\alpha_0,\beta_0}} \frac{\vert
t\alpha_0-\beta_0-x\vert}{\sqrt{1+t^2-(\beta_0-t\alpha_0+x)^2}}\
dt \, dx =\int\limits_{\DD_{0,0}} \frac{\vert
v\vert}{\sqrt{1+u^2-v^2}} \ du\, dv \\
& =\sqrt{2}+\ln (1+\sqrt{2}\, ).
\end{split}
\end{equation*}

(ii) The same change of variable as in (i) gives
\begin{equation*}
\iint\limits_{\DD_{\alpha_0,\beta_0}} \sqrt{\Phi (t,x)}\ dt\, dx=
\iint\limits_{\DD_{0,0}} \sqrt{1+u^2-v^2}\ du\, dv=\frac{2\pi}{3}
\ . \end{equation*} \qed

We start to evaluate $S_Q$. If $D\in J_{q,a}$, then $D\in \Omega_q
=\big( -(2+r_0\sqrt{2})q,(2+r_0\sqrt{2})q\big)$. The equality
$aq^\prime -a^\prime q=D$ is equivalent to $a^\prime
=\frac{aq^\prime -D}{q}$. Hence in the inner sum we should sum
over $q^\prime \in [q,Q]$ such that $aq^\prime =D\pmod{q}$ and
\begin{equation*}
a\leq \frac{aq^\prime -D}{q}\leq Q,
\end{equation*}
or equivalently
\begin{equation}\label{3.5}
\max\bigg\{ q,q+\frac{D}{a}\bigg\} \leq q^\prime \leq \min \bigg\{
Q,\frac{qQ+D}{a}\bigg\} .
\end{equation}

Next we show that the bulk of the contribution to $S_Q$ of
\eqref{3.5} only comes from
\begin{equation*}
q\leq q^\prime \leq Q.
\end{equation*}
To see this, notice first that for $q^\prime$ fixed, the relations
$D=aq^\prime \hspace{-4pt}\pmod{q}$ and $\vert D\vert \leq
(r_0+\sqrt{2})q$ imply that $D$ takes at most $3+[r_0]$ values. So
the total contribution of terms with $0\leq q-a\leq 2+[r_0]$ to
$S_Q$ is
\begin{equation*}
\ll \frac{1}{Q^4} \sum\limits_{q=1}^Q \sum\limits_{q^\prime =1}^Q
(3+[r_0])\cdot \frac{(q^\prime-q)q^\prime}{q^2}\cdot \sqrt{2} r_0
q\ll_{r_0} \frac{1}{Q} \sum\limits_{q=1}^Q \frac{1}{q} \ll
\frac{\ln Q}{Q}\ .
\end{equation*}
When $q-a\geq 3+[r_0]$, we get $-D\leq (3+[r_0])Q\leq (q-a)Q$,
thus $Q\leq \frac{qQ+D}{a}$. Suppose now that $q^\prime$ is
between $q$ and $q+\frac{D}{a}$. Owing again to $aq^\prime
=D\pmod{q}$, it follows that $aq^\prime =D+kq$ for some integer
$k$. But the range of $q^\prime$ is an interval of length
$\frac{\vert D\vert}{a}$, hence the range of
$k=\frac{aq^\prime}{q}-\frac{D}{q}$ is an interval of length $\ll
1+\frac{a}{q}\cdot \frac{\vert D\vert}{a}\leq 4+[r_0]$. Thus $k$,
and consequently $q^\prime$, take at most $O(1)$ values. Besides,
in this case we have $0\leq q^\prime -q\leq \frac{\vert
D\vert}{a}$. Hence the contribution to $S_Q$ of terms with
$q^\prime$ between $q$ and $q+\frac{D}{a}$ is
\begin{equation*}
\ll \frac{1}{Q^4} \sum\limits_{1\leq a\leq q\leq Q}
\sum\limits_{D\in \Omega_q} \frac{\frac{\vert D\vert}{a}\cdot
q}{q^2} \cdot q =\frac{1}{Q^4} \sum\limits_{1\leq a\leq q\leq Q}
\frac{1}{a} \sum\limits_{D\in \Omega_q} \vert D\vert \ll
\frac{1}{Q^4} \sum\limits_{a=1}^Q \frac{1}{a} \sum\limits_{q=1}^Q
q^2 \ll \frac{\ln Q}{Q} \ .
\end{equation*}
Therefore we have shown that we can replace the summation
conditions in the inner sum from \eqref{3.2} by $q^\prime \in [q,
Q]$ and $aq^\prime =D\hspace{-4pt}\pmod{q}$.

We write $x_0=r_0 \alpha_0$ and $y_0=r_0\beta_0$. Take
$\DD=\DD_{\alpha_0,\beta_0}$, $\Phi=\Phi_{\alpha_0,\beta_0}$ and
$\psi=\psi_{\alpha_0,\beta_0}$, unless otherwise specified, and
note that
\begin{equation*}
r_0^2 (q^2+a^2)-(y_0q-x_0 a+D)^2 =r_0^2 q^2 \Phi \bigg(
\frac{a}{q},\frac{D}{r_0 q}\bigg) .
\end{equation*}
Then \eqref{3.1} and the previous considerations lead to
\begin{equation*}
\begin{split}
S_Q & =\frac{2r_0}{Q^4} \sum\limits_{0\leq a\leq q\leq Q}\
\sum\limits_{D\in r_0 qJ_{\frac{a}{q}}}\
\sum\limits_{\substack{q^\prime \in [q, Q]
\\ aq^\prime =D \hspace{-8pt}\pmod{q}}} \frac{(q^\prime
-q)q^\prime}{q}\ \sqrt{\Phi \bigg(
\frac{a}{q},\frac{D}{r_0q}\bigg)}\  +O(Q^{-1}\ln Q) \\
& =\frac{2r_0}{Q^4} \sum\limits_{q=1}^Q \frac{1}{q}
\sum\limits_{D\in r_0 q\mathrm{pr}_2
\DD}\sum\limits_{\substack{q^\prime \in [q,Q]
\\ a\in qI_{\hspace{-2pt}\frac{D}{r_0q}} \\ aq^\prime
=D\hspace{-8pt}\pmod{q}}} (q^\prime -q)q^\prime \sqrt{\Phi
\bigg(\frac{a}{q},\frac{D}{r_0 q}\bigg)} \ +O(Q^{-1}\ln Q),
\end{split}
\end{equation*}
where $I_x$ and $J_t$ are defined as in \eqref{3.3} and
\eqref{3.4}.

Take $d=\gcd (q,q^\prime)$ and write $q=dq_0$, $q^\prime
=dq_0^\prime$. Then $d$ divides $D$, so $D=dD_0$. The congruence
$adq_0^\prime =D\pmod{dq_0}$ is equivalent to $aq_0^\prime
=D_0\pmod{q_0}$, and we may write
\begin{equation}\label{3.6}
S_Q=\frac{2r_0}{Q^4} \sum\limits_d \sum\limits_{q_0=1}^{\big[
\frac{Q}{d}\big]} \frac{d}{q_0} \sum\limits_{D_0\in r_0
q_0\mathrm{pr}_2 \DD} \hspace{-10pt}
\sum\limits_{\substack{q_0^\prime \in \big[
q_0,\frac{Q}{d}\big] \\ \gcd (q_0^\prime,q_0)=1 \\
a\in dq_0 I_{\hspace{-3pt} \frac{D_0}{r_0q_0}} \\
aq_0^\prime =D_0 \hspace{-8pt} \pmod{q_0}}} \hspace{-15pt}
(q_0^\prime -q_0)q_0^\prime \sqrt{\Phi \bigg(
\frac{a}{dq_0},\frac{D_0}{r_0q_0}\bigg)} \ +O(Q^{-1}\ln Q).
\end{equation}

\medskip
To estimate the inner sum above, we need some information about
the distribution of solutions of the congruence $xy=h\hspace{-5pt}
\pmod{q}$. We shall employ the following result, which is a
consequence of Proposition \ref{PA2} from the Appendix.

\medskip

\begin{prop}\label{L3.3}
Assume that $q\geq 1$ and $h$ are two given integers, $\II$ and
$\JJ$ are intervals, and $f:\II \times \JJ \rightarrow \R$ is a
$C^1$ function. Then for every integer $T>1$ and every $\delta >0$
\begin{equation*}
\begin{split}
\sum_{\substack{a\in \II,\, b\in \JJ \\
ab=h\hspace{-8pt}\pmod{q}\\ \gcd (b,q)=1}} \hspace{-3pt} & f(a,b)
=\frac{\varphi (q)}{q^2} \iint\limits_{\II \times \JJ} f(x,y) dx
dy +O_\delta \Bigg( \Big( 1+\frac{\vert \II\vert}{q}\Big) \Big(
1+\frac{\vert \JJ\vert}{q}\Big) T^2 \| f\|_\infty
q^{\frac{1}{2}+\delta}\gcd (h,q)   \\
& \qquad \qquad + \Big( 1+\frac{\vert \II\vert}{q}\Big) \Big(
1+\frac{\vert \JJ\vert}{q}\Big) T\| Df\|_\infty
q^{\frac{3}{2}+\delta} \gcd(h,q)+ \frac{\vert \II\vert\, \vert
\JJ\vert \| Df\|_\infty}{T} \Bigg),
\end{split}
\end{equation*}
where $\| \cdot \|_\infty =\| \cdot \|_{\infty, \II \times \JJ}$.
\end{prop}

We now return to the formula for $S_Q$ given in \eqref{3.6} and
first give an upper bound for the contribution to $S_Q$ of
quadruples $(d,q_0,D_0,a)$ for which
\begin{equation}\label{3.7}
0\leq r_0^2d^2q_0^2 \Phi \bigg(
\frac{a}{dq_0},\frac{D_0}{r_0q_0}\bigg)
=r_0^2(a^2+d^2q_0^2)-(dq_0y_0-ax_0+dD_0)^2 \leq L^2,
\end{equation}
with $L=L_{q_0}>1$ to be chosen later.

\medskip

\begin{lem}\label{L3.4}
Let $F(a)=ua^2+va+w$ with $u\neq 0$. Then for any $K$ and $L$
\begin{equation*}
\big| \{ a\in \R \,:\,K\leq F(a)\leq K+L^2 \}\big| \leq
\frac{2\vert L\vert}{\sqrt{\vert u\vert}} \ .
\end{equation*}
\end{lem}

{\sl Proof.} Using
\begin{equation*}
\{ a\,:\,K\leq F(a)\leq K+L^2\} =\{ a\,:\, -K-L^2\leq -F(a)\leq
-K\}
\end{equation*}
we see that it suffices to consider the case $u>0$. In this case
the statement follows from the fact that the double inequality
$K\leq F(t)\leq K+L^2$ is equivalent to
\begin{equation*}
\frac{1}{u} \bigg( K+\frac{v^2-4uw}{4u}\bigg) \leq \bigg(
a+\frac{v}{2u}\bigg)^2 \leq \frac{L^2}{u} +\frac{1}{u} \bigg(
K+\frac{v^2-4uw}{4u}\bigg) ,
\end{equation*}
and from the inequality
\begin{equation*}
\vert \sqrt{x}-\sqrt{y} \, \vert \leq \sqrt{\vert x-y\vert} \ .
\end{equation*}
\qed

Suppose that $(d,q_0,D_0)$ is fixed and consider the following two
cases:

{\bf Case 1)} $\ r_0\neq x_0$.

By \eqref{3.7} and Lemma \ref{L3.4} the range of $a$ is the union
of at most two intervals of length $\frac{L_{q_0}}{r_0}$. Hence
$a$ can only assume $O(L_{q_0})$ values. But, for each $a$,
$q_0^\prime$ belongs to $\big[ q_0,\frac{Q}{d}\big]$ and is
subject to the condition $q_0^\prime a=D_0 \hspace{-4pt}
\pmod{q_0}$. Hence $q_0^\prime$ takes $O\big(
1+\frac{Q}{dq_0}\big)=O\big( \frac{Q}{dq_0}\big)$ values. Thus the
contribution to $S_Q$ of quadruples $(d,q_0,D_0,a)$ which satisfy
\eqref{3.7} is
\begin{equation}\label{3.8}
\ll \frac{1}{Q^4} \sum\limits_d \sum\limits_{q_0} \frac{d}{q_0}
\sum\limits_{\vert D_0\vert \ll q_0} L_{q_0} \cdot \frac{Q}{dq_0}
\bigg( \frac{Q}{d}\bigg)^2 \frac{L_{q_0}}{dq_0} \ll \frac{1}{Q}
\sum\limits_{q_0} \frac{L_{q_0}^2}{q_0^2} \ .
\end{equation}

{\bf Case 2)} $\ r_0=x_0$ thus $\alpha_0=1$.

In this case we collect directly from \eqref{3.7}
\begin{equation*}
\frac{(dq_0y_0+dD_0)^2-r_0^2d^2q_0^2}{2r_0d(D_0+q_0y_0)} \leq
a\leq
\frac{L^2+(dq_0y_0+dD_0)^2-r_0^2d^2q_0^2}{2r_0d(D_0+q_0y_0)}\ .
\end{equation*}
Hence $a$ can only assume $O\big( \frac{L_{q_0}^2}{dD_0}\big)$
values and we find, arguing as in Case 1, that the contribution to
$S_Q$ of quadruples $(d,q_0,D_0,a)$ which satisfy \eqref{3.7} is
\begin{equation}\label{3.9}
\begin{split}
\ll \frac{1}{Q^4} \sum\limits_d \sum\limits_{q_0} \frac{d}{q_0}
\sum\limits_{\vert D_0\vert \ll q_0} \frac{L_{q_0}^2}{dD_0} \cdot
\frac{Q}{dq_0}\bigg( \frac{Q}{d}\bigg)^2 \frac{L_{q_0}}{dq_0} &
\ll \frac{1}{Q} \sum\limits_{q_0=1}^Q \frac{\log q_0}{q_0^3} \cdot
L_{q_0}^3 \\ & \ll_\delta \frac{1}{Q} \sum\limits_{q_0=1}^Q
\frac{L_{q_0}^3}{q_0^{3-\delta}} \ .
\end{split}
\end{equation}

Next we investigate the case
\begin{equation*}
r_0^2d^2q_0^2 \Phi \bigg( \frac{a}{dq_0},\frac{D_0}{r_0q_0}\bigg)
=r_0^2(a^2+d^2q_0^2)-(dq_0y_0-ax_0+dD_0)^2 \geq L^2.
\end{equation*}
We consider the range of $q_0^\prime$:
\begin{equation*}
\II_{q_0,d}=\bigg[ q_0,\frac{Q}{d}\bigg] ,
\end{equation*}
the range of $a$ (which is the union of at most two intervals):
\begin{equation*}
\JJ_{q_0,D_0,d,L}=\bigg\{ y\in [0,dq_0]\, :\, \Phi \bigg(
\frac{y}{dq_0},\frac{D_0}{r_0q_0}\bigg) \geq \frac{L^2}{r_0^2d^2
q_0^2}\bigg\} \subseteq dq_0 I_{\hspace{-3pt}\frac{D_0}{r_0q_0}} \
,
\end{equation*}
and the functions
\begin{equation*}
\begin{split}
G(x)=G_{q_0}(x) & :=(x-q_0)x,\qquad x\in \II_{q_0,d} ,\\
\Psi(y)=\Psi_{q_0,D_0,d}(y) & :=\sqrt{\Phi \bigg(
\frac{y}{dq_0},\frac{D_0}{r_0q_0}\bigg)}\ ,
\qquad y\in \JJ_{q_0,D_0,d,L} \subseteq [0,dq_0] ,\\
F(x,y)=F_{q_0,D_0,d}(x,y) & :=G(x) \Psi(y), \qquad (x,y)\in
\II_{q_0,d} \times \JJ_{q_0,D_0,d,L} .
\end{split}
\end{equation*}
With this notation the following estimates hold on
$\II_{q_0,d}\times \JJ_{q_0,D_0,d,L}$:
\begin{equation*}
\begin{split}
& \| G\|_\infty \ll \frac{Q^2}{d^2} \, ,\qquad \| \Psi \|_\infty
\ll 1 ,\quad \| F\|_\infty \leq \| G\|_\infty \|
\Psi \|_\infty \ll \frac{Q^2}{d^2} \, , \\
& \| G^\prime \|_\infty \ll \frac{Q}{d} \, , \qquad \| \Psi^\prime
\|_\infty =\sup\limits_{y\in \JJ_{q_0,D_0,d,L}} \frac{\big|
\frac{(1-\alpha_0^2)y}{d^2q_0^2} + \frac{\alpha_0}{dq_0} \big(
\frac{D_0}{r_0q_0}+\beta_0\big)\big|}{\sqrt{\Phi \big(
\frac{y}{dq_0},\frac{D_0}{r_0q_0}\big)}} \ \ll
\frac{\frac{1}{dq_0}}{\frac{L_{q_0}}{dq_0}} =\frac{1}{L_{q_0}} \, ,\\
 & \| DF\|_\infty \leq \| G\|_\infty \|
\Psi^\prime \|_\infty +\| G^\prime \|_\infty \| \Psi \|_\infty \ll
\frac{Q^2}{d^2}\cdot \frac{1}{L_{q_0}} +\frac{Q}{d}
=\frac{Q^2}{dL_{q_0}} \bigg( \frac{1}{d}+\frac{L_{q_0}}{Q}\bigg)
\\ & \hspace{7cm}
\leq \frac{Q^2}{dL_{q_0}} \bigg( \frac{1}{d}+\frac{q_0}{Q}\bigg)
\ll \frac{Q^2}{d^2 L_{q_0}} \ .
\end{split}
\end{equation*}

Applying Proposition \ref{L3.3} we find that
\begin{equation*}
\begin{split}
\sum\limits_{\substack{q_0^\prime \in \big[ q_0,\frac{Q}{d}\big]
\\ \gcd (q_0^\prime,q_0)=1
\\ a\in \JJ_{q_0,D_0,d,L} \\ aq_0^\prime =D_0 \hspace{-8pt}\pmod{q_0}}}
\hspace{-15pt} (q_0^\prime-q_0)q_0 & \sqrt{\Phi \bigg(
\frac{a}{dq_0},\frac{D_0}{r_0q_0}\bigg)} =
\sum\limits_{\substack{q_0^\prime \in \big[ q_0,\frac{Q}{d}\big]
\\ \gcd (q_0^\prime,q_0)=1
\\ a\in \JJ_{q_0,D_0,d,L} \\ aq_0^\prime =D_0 \hspace{-8pt}\pmod{q_0}}}
\hspace{-15pt} G_{q_0}(q_0^\prime) \Psi_{q_0,D_0,d}(a)\\ &
=\frac{\varphi(q_0)}{q_0^2} \int\limits_{q_0}^{\frac{Q}{d}}
G_{q_0}(x)\, dx \hspace{-10pt} \int\limits_{\JJ_{q_0,D_0,d,L}}
\hspace{-15pt} \Psi_{q_0,D_0,d}(y)\, dy +\EE_{q_0,D_0,d,L},
\end{split}
\end{equation*}
where
\begin{equation*}
\begin{split}
\EE_{q_0,D_0,d,L} & \ll_\delta \frac{Q}{dq_0}\cdot dT^2
q_0^{\frac{1}{2}+\delta} \cdot \frac{Q^2}{d^2}\cdot \gcd
(D_0,q_0)+\frac{Q}{dq_0}\cdot dTq_0^{\frac{3}{2}+\delta} \cdot
\frac{Q^2}{d^2L}\gcd (D_0,q_0) \\
& \qquad \quad +\frac{Q}{d}\cdot dq_0\cdot \frac{Q^2}{d^2 L}\cdot
\frac{1}{T} \\
& =\frac{Q^3}{d^2 q_0} \bigg( T^2 q_0^{\frac{1}{2}+\delta} \gcd
(D_0,q_0)+\frac{Tq_0^{\frac{3}{2}+\delta} \gcd
(D_0,q_0)}{L}+\frac{q_0^2}{LT}\bigg).
\end{split}
\end{equation*}
Using
\begin{equation*}
\sum\limits_{\vert D_0\vert \ll q_0} \gcd (D_0,q_0)\ll_\delta
q_0^{1+\delta},
\end{equation*}
we see that
\begin{equation*}
\sum\limits_{\vert D_0 \vert \ll q_0} \EE_{q_0,D_0,d,L} \ll_\delta
\frac{Q^3}{d^2 q_0} \bigg( T^2
q_0^{\frac{3}{2}+\delta}+\frac{Tq_0^{\frac{5}{2}+\delta}}{L}+
\frac{q_0^3}{LT} \bigg).
\end{equation*}
Taking $T=[q_0^{\frac{1}{5}}]$ and $L=q_0^{\frac{9}{10}}$ we find
that
\begin{equation*}
\sum\limits_{\vert D_0 \vert \ll q_0} \EE_{q_0,D_0,d,L} \ll_\delta
\frac{Q^3}{d^2 q_0} \cdot q_0^{2-\frac{1}{10}+\delta}
=\frac{Q^3}{d^2}\cdot q_0^{1-\frac{1}{10}+\delta} .
\end{equation*}
Thus the total contribution of $\EE_{q_0,D_0,d,L}$ to $S_Q$ is
\begin{equation*}
\ll_\delta \frac{1}{Q} \sum\limits_{d=1}^Q \sum\limits_{q_0=1}^Q
\frac{d}{q_0}\cdot \frac{1}{d^2}\cdot q_0^{1-\frac{1}{10}+\delta}
\ll Q^{-\frac{1}{10}+2\delta}.
\end{equation*}
Moreover, the quantities in \eqref{3.8} and \eqref{3.9} become
\begin{equation*}
\ll \frac{1}{Q} \sum\limits_{q_0=1}^Q q_0^{-\frac{1}{5}} \ll
\frac{Q^{\frac{4}{5}}}{Q}=Q^{-\frac{1}{5}}\leq Q^{-\frac{1}{10}}\
,
\end{equation*}
and respectively
\begin{equation*}
\ll \frac{1}{Q} \sum\limits_{q_0=1}^Q
q_0^{-\frac{3}{10}+\delta}\ll Q^{-\frac{3}{10}+\delta} \leq
Q^{-\frac{1}{10}+\delta} .
\end{equation*}
Thus we gather
\begin{equation}\label{3.10}
S_Q=M_Q+O_\delta(Q^{-\frac{1}{10}+\delta}),
\end{equation}
with
\begin{equation*}
M_Q=\frac{2r_0}{Q^4} \sum\limits_d \sum\limits_{q_0=1}^{\big[
\frac{Q}{d}\big]} \frac{d}{q_0} \sum\limits_{D_0 \in r_0 q_0
\mathrm{pr}_2 \DD} \frac{\varphi(q_0)}{q_0^2}
\int\limits_{q_0}^{\frac{Q}{d}} G_{q_0}(x)\, dx \hspace{-15pt}
\int\limits_{\JJ_{q_0,D_0,d,q_0^{9/10}}} \hspace{-20pt}
\Psi_{q_0,D_0,d}(y)\  dy.
\end{equation*}

Next we show that we can replace $\JJ_{q_0,D_0,d,q_0^{9/10}}$ by
$dq_0 I_{\hspace{-3pt} \frac{D_0}{r_0q_0}}$ in the last integral.
Clearly $\JJ_{q_0,D_0,d,q_0^{9/10}} \subseteq dq_0
I_{\hspace{-3pt}\frac{D_0}{r_0q_0}}$ and by Lemma \ref{L3.4} we
have
\begin{equation*}
\Big|\, dq_0 I_{\hspace{-3pt}\frac{D_0}{r_0q_0}} \setminus
\JJ_{q_0,D_0,d,q_0^{9/10}}\Big| \ll
2\sqrt{\frac{(q_0^{9/10})^2}{r_0^2-x_0^2}} \ \ll
q_0^{\frac{9}{10}} .
\end{equation*}
Thus
\begin{equation*}
0\leq \int\limits_{dq_0I_{\hspace{-3pt}\frac{D_0}{r_0q_0}}
\setminus \JJ_{q_0,D_0,d,q_0^{9/10}}} \hspace{-30pt}
\Psi_{q_0,D_0,d}(y)\, dy  \leq  \Big|\, dq_0
I_{\hspace{-3pt}\frac{D_0}{r_0q_0}} \setminus
\JJ_{q_0,D_0,d,q_0^{9/10}}\Big| \cdot \frac{q_0^{9/10}}{dq_0} \ll
\frac{q_0^{4/5}}{d}\ ,
\end{equation*}
and as a result the error that results from replacing
$\JJ_{q_0,D_0,d,L}$ by $dq_0I_{\hspace{-3pt}\frac{D_0}{r_0q_0}}$
in $M_Q$ is
\begin{equation*}
\ll \frac{1}{Q^4} \sum\limits_d \sum\limits_{q_0=1}^{\big[
\frac{Q}{d}\big]} \frac{d\varphi(q_0)}{q_0^3} \cdot
\frac{q_0^{4/5}}{d} \int\limits_{q_0}^{\frac{Q}{d}} G_{q_0}(x)\,
dx.
\end{equation*}
On the other hand we find that
\begin{equation*}
\int\limits_{q_0}^{\frac{Q}{d}} G_{q_0}(x)\, dx
=\int\limits_{0}^{\frac{Q}{d}-q_0} x(x+q_0)\, dx=q_0^3 \  g\bigg(
\frac{Q}{dq_0} \bigg),
\end{equation*}
where
\begin{equation*}
g(t)=\frac{(t-1)^3}{3}+\frac{(t-1)^2}{2} \ .
\end{equation*}
In particular, the integral of $G_{q_0}$ on $\big[
q_0,\frac{Q}{d}\big]$ is $\ll \frac{Q^3}{d^3}$, and we find that
the total cost of replacing $\JJ_{q_0,D_0,d,q_0^{9/10}}$ by $dq_0
I_{\hspace{-3pt}\frac{D_0}{r_0q_0}}$ in $M_Q$ is
\begin{equation*}
\ll \frac{1}{Q^4} \sum\limits_d \sum\limits_{q_0=1}^{\big[
\frac{Q}{d}\big]} \frac{\varphi(q_0)}{q_0^3}\cdot
q_0^{\frac{4}{5}} \cdot \frac{Q^3}{d^3} \leq \frac{1}{Q}
\sum\limits_{d=1}^\infty \frac{1}{d^3} \sum\limits_{q_0=1}^\infty
q_0^{-\frac{6}{5}} \ll \frac{1}{Q} \  .
\end{equation*}
Thus we infer that
\begin{equation}\label{3.11}
\begin{split}
M_Q & =\frac{2r_0}{Q^4} \sum\limits_d \sum\limits_{q_0=1}^{\big[
\frac{Q}{d}\big]} \frac{d}{q_0}\cdot\frac{\varphi(q_0)}{q_0^2}
\cdot q_0^3 \, g\bigg( \frac{Q}{dq_0}\bigg)
\hspace{-2pt}\sum\limits_{D_0 \in r_0q_0 \mathrm{pr}_2 \DD}\
\int\limits_{dq_0I_{\hspace{-3pt}\frac{D_0}{r_0q_0}}}
\hspace{-10pt} \Psi_{q_0,D_0,d}(y)\, dy +O(Q^{-1}) \\
 & =\frac{2r_0}{Q^4} \sum\limits_d \sum\limits_{q_0=1}^{\big[
\frac{Q}{d}\big]} d\varphi(q_0)\, g\bigg( \frac{Q}{dq_0}\bigg)
\sum\limits_{D_0 \in r_0q_0 \mathrm{pr}_2 \DD}\ \
\int\limits_{dq_0I_{\hspace{-3pt}\frac{D_0}{r_0q_0}}}
\hspace{-10pt} \sqrt{\Phi \bigg(
\frac{y}{dq_0},\frac{D_0}{r_0q_0}\bigg)}\  dy +O(Q^{-1}) \\
& =\frac{2r_0}{Q^4} \sum\limits_d \sum\limits_{q_0=1}^{\big[
\frac{Q}{d}\big]} d^2 q_0 \varphi(q_0) \, g\bigg(
\frac{Q}{dq_0}\bigg) \sum\limits_{D_0 \in r_0q_0 \mathrm{pr}_2
\DD}\ \int\limits_{I_{\hspace{-3pt}\frac{D_0}{r_0q_0}}}
\hspace{-7pt}
\sqrt{\Phi \bigg( t,\frac{D_0}{r_0 q_0}\bigg)} \  dt +O(Q^{-1}) \\
& =\frac{2r_0}{Q^4} \sum\limits_d \sum\limits_{q_0=1}^{\big[
\frac{Q}{d}\big]} d^2 q_0 \varphi(q_0) \, g\bigg(
\frac{Q}{dq_0}\bigg) \sum\limits_{D_0 \in r_0q_0 \mathrm{pr}_2
\DD} \hspace{-8pt} \psi \bigg( \frac{D_0}{r_0q_0}\bigg)
+O(Q^{-1}).
 \end{split}
\end{equation}

By Euler-MacLaurin summation and Lemma \ref{L3.0} the inner sum
above is given by
\begin{equation*}
\begin{split}
& \int\limits_{r_0q_0\mathrm{pr}_2 \DD} \psi \bigg(
\frac{u}{r_0q_0}\bigg) du+O\Bigg( \sup\limits_{x\in \mathrm{pr}_2
\DD} \vert \psi(x)\vert +\int\limits_{r_0q_0\mathrm{pr}_2 \DD}
\bigg| \frac{1}{r_0q_0}\cdot
\psi^\prime \bigg( \frac{u}{r_0q_0}\bigg) \bigg|\, du\Bigg) \\
& \qquad =r_0q_0 \int\limits_{\mathrm{pr}_2 \DD} \psi (v)dv+O(1)
=\frac{2\pi r_0q_0}{3}+O(1).
\end{split}
\end{equation*}

Inserting this back into \eqref{3.11} we obtain
\begin{equation}\label{3.12}
\begin{split}
M_Q & =\frac{4\pi r_0^2}{3Q^4} \sum\limits_d
\sum\limits_{q_0=1}^{\big[ \frac{Q}{d}\big]} \varphi(q_0) (dq_0)^2
g\bigg( \frac{Q}{dq_0}\bigg) +O\bigg( \frac{\ln Q}{Q}\bigg) \\
\big( q=dq_0\in [1,Q]\big) \qquad & =\frac{4\pi r_0^2}{3Q^4}
\sum\limits_{q=1}^Q q^2 g\bigg( \frac{Q}{q}\bigg) \sum\limits_{q_0
\vert q} \varphi(q_0) +O\bigg( \frac{\ln^2 Q}{Q}\bigg) \\
& =\frac{4\pi r_0^2}{3Q^4} \sum\limits_{q=1}^Q q^3 g\bigg(
\frac{Q}{q}\bigg) +O\bigg( \frac{\ln^2 Q}{Q}\bigg) \\
& =\frac{4\pi r_0^2}{3} \int\limits_0^1 \bigg(
\frac{(1-x)^3}{3}+\frac{x(1-x)^2}{2} \bigg) dx+O\bigg(
\frac{\ln^2 Q}{Q} \bigg) \\
& =\frac{\pi r_0^2}{6}+O\bigg( \frac{\ln^2 Q}{Q}\bigg) .
\end{split}
\end{equation}
Proposition \ref{P3.1} now follows from \eqref{3.10} and
\eqref{3.12}.

\bigskip

\section{On the $6$-level correlations}
In this section we prove Theorem \ref{T1.2}. We first prove a
counting result.

\medskip

\begin{lem}\label{L4.1}
Let $a,b,d$ be positive integers with $\gcd (a,b,d)=1$. Let $d_1$
denote the largest divisor of $d$ which is relatively prime with
$b$, and put $d_2=\frac{d}{d_1}$. Then
\begin{equation}\label{4.1}
\# \{ 0\leq m\leq 2d-1\, :\, \gcd (a+bm,d)=1\} =2\varphi(d_1)d_2.
\end{equation}
\end{lem}

{\sl Proof.} Using M\" obius inversion we express the left-hand
side of \eqref{4.1} as
\begin{equation*}
\sum\limits_{0\leq m\leq 2d-1} \sum\limits_{\substack{D\vert d \\
D\vert a+bm}} \mu (D)=\sum\limits_{D\vert d} \mu (D)\ \# \{ 0\leq
m\leq 2d-1\, :\, D\vert \, a+bm\} .
\end{equation*}

Note that if $D$ does not divide $d_1$, then there is no $m$ for
which $D\vert\, a+bm$. Indeed, if for some $m$ we have $D\vert\,
a+bm$, then $a=Dk-bm$ for some $k\in \Z$. Hence $a$ is divisible
by $\gcd (D,b)$. Then $\gcd (D,b)$ divides $\gcd (a,b,d)=1$, so
$\gcd (D,b)=1$, and by the definition of $d_1$ it follows that $D$
divides $d_1$. Therefore
\begin{equation*}
\# \{ 0\leq m\leq 2d-1 \, :\, \gcd (a+bm,d)=1\}
=\sum\limits_{D\vert\, d_1} \mu(D)\ \# \{ 0\leq m\leq 2d-1\, :\,
D\vert \, a+bm\}.
\end{equation*}

For $D\vert\, d_1$ we have $\gcd (D,b)=1$ and there is a unique
solution $m\pmod{D}$ to the congruence $a+bm=0\pmod{D}$. Thus
there are exactly $\frac{2d}{D}$ values of $m$ in $\{
0,1,\dots,2d-1\}$ for which $D\vert\, a+bm$. As a result we infer
that
\begin{equation*}
\begin{split}
& \# \{ 0\leq m\leq 2d-1\, :\, \gcd (a+bm,d)=1\}
=\sum\limits_{D\vert\, d_1} \mu(D)\
\frac{2d}{D}=2d\sum\limits_{D\vert\, d_1} \frac{\mu(D)}{D} \\
& \qquad \qquad =2d\ \frac{\varphi(d_1)}{d_1} =2\varphi(d_1)d_2 ,
\end{split}
\end{equation*}
which proves the lemma. \qed

\bigskip

For any positive integers $a,b,q$, consider the set
\begin{equation*}
\NN_{a,b,q}=\big\{ (A,B)\, :\ 1\leq A,B\leq 2q,\ \gcd (A,B)=1, \
q\vert\, Ab-Ba\big\}.
\end{equation*}

\medskip

\begin{lem}\label{L4.2}
For $q$ large and $1\leq a,b\leq q$ such that $\gcd (a,b,q)=1$
\begin{equation*}
\# \NN_{a,b,q} \gg \frac{\varphi(q)}{\ln q} \gg \frac{q}{\ln q \ln
\ln q}\ .
\end{equation*}
\end{lem}

\medskip

{\bf Remark.} If in the definition of $\NN_{a,b,q}$ we took the
range of $A$ and $B$ to be $[1,q]$ instead of $[1,2q]$, the
cardinality of $\NN_{a,b,q}$ would be much smaller. For example,
if $1\leq a\leq q$, $\gcd (a,q)=1$, and $b=a$, then $q\vert\,
A-B$, and in the range $1\leq A,B\leq q$ this forces $A=B$. Then
the only pair $(A,B)$ with $\gcd (A,B)=1$ is $(1,1)$, so
$\NN_{a,b,q}$ will only contain one element. Lemma \ref{L4.2}
shows a sudden increase in the cardinality of $\NN_{a,b,q}$ when
the range of $A$ and $B$ increases by a factor $2$.

\medskip

{\sl Proof of Lemma \ref{L4.2}.}  To make a choice, assume in what
follows next that $d=\gcd (a,q)\leq \gcd (b,q)$. Then $\gcd
(d,b)=1$ since $\gcd (q,a,b)=1$. It follows that for any solution
$(A,B)$ to the congruence $Ab=Ba\pmod{q}$, $A$ has to be divisible
by $d$. Write $q=dq_1$ and $a=da_1$, so $\gcd (a_1,q_1)=1$. Denote
by $\bar{a}_1$ the multiplicative inverse of $a_1\pmod{q_1}$ in
the interval $[1,q_1]$.

Note that since $q$ is divisible by the product $\gcd (a,q)\gcd
(b,q)$, we have $d<\sqrt{q}$. Therefore $q_1>\sqrt{q}$. So $q_1$
is large for large $q$, and by the Prime Number Theorem we know
that
\begin{equation}\label{4.2}
\# \{ p\ \mbox{\rm prime}\, :\, q_1 <p\leq 2q_1\} \sim
\frac{q_1}{\ln q_1} \ .
\end{equation}

For any fixed prime $p$ with $q_1<p\leq 2q_1$, we count the
solutions $1\leq B\leq 2q$ of the congruence
\begin{equation}\label{4.3}
dpb=Ba\pmod{q}.
\end{equation}
This is equivalent to
\begin{equation}\label{4.4}
pb=Ba_1 \pmod{q_1}.
\end{equation}
Since $\gcd (q_1,a_1)=1$, this congruence has a unique solution
modulo $q_1$, namely $B=\overline{a}_1 pb\pmod{q_1}$. Denote by
$B_0$ the solution to \eqref{4.4} which belongs to the interval
$[1,\leq q_1]$. Then \eqref{4.3} will have $2d$ solutions, given
by
\begin{equation}\label{4.5}
B=B_0+q_1 m,\qquad 0\leq m\leq 2d-1.
\end{equation}

It remains to be seen how many of the numbers $B$ from \eqref{4.5}
are relatively prime with $A=dp$. Note first that at most two such
numbers $B$ are divisible by $p$. Assume that
\begin{equation}\label{4.6}
p\vert\, B_0+q_1 m_1,\quad p\vert\, B_0+q_1 m_2,\quad \mbox{\rm
and}\quad p\vert \, B_0+q_1m_3,
\end{equation}
with $0\leq m_1<m_2<m_3\leq 2d-1$. Since $p>q_1$, $p$ does not
divide $q_1$. Then it follows from \eqref{4.6} that
\begin{equation}\label{4.7}
p\vert\, m_2-m_1 \qquad \mbox{\rm and}\qquad p\vert\, m_3-m_2 .
\end{equation}
Here at least one of the differences $m_2-m_1$, $m_3-m_2$ is less
than $d$, and $d<\sqrt{q}<q_1<p$, so it is impossible that both
divisibilities in \eqref{4.7} hold true. Therefore at most two
numbers from \eqref{4.5} are divisible by $p$. Note also, by the
same reasoning, that for smaller values of $d$ - more precisely
for $d<\sqrt{\frac{q}{2}}$, that is $d<\frac{q_1}{2}$ - one has
$2d<p$. Then one concludes that at most one number $B$ from
\eqref{4.5} can be a multiple of $p$.

We now count the numbers $B$ from \eqref{4.5} which are relatively
prime with $d$. We claim that $\gcd(B_0,q_1,d)=1$. Indeed, let us
assume this fails and choose a prime divisor $p_1$ of $\gcd
(B_0,q_1,d)$. Since $p_1\vert\, d$, we have $p_1\vert\, a$ and
$p_1\vert\, q$. Recall that $B_0$ satisfies \eqref{4.4}, hence
\begin{equation}\label{4.8}
pb=B_0a_1+q_1 k
\end{equation}
for some $k\in \Z$. Here $p_1\vert\, B_0$, $p_1\vert\, q_1$, so
$p_1$ must also divide the left side of \eqref{4.8}. The
inequalities $p_1\leq B_0\leq q_1<p$ show that $p_1\neq p$, so $p$
divides $b$. But then $p_1\vert\, \gcd(a,b,q)=1$, and we obtain a
contradiction. This shows that $\gcd (B_0,q_1,d)=1$. Then Lemma
\ref{L4.1} is applicable to $B_0,q_1,d$. If $d_1$ denotes the
largest divisor of $d$ which is relatively prime with $q_1$, and
$d_2=\frac{d}{d_1}$, then Lemma \ref{L4.1} provides
\begin{equation*}
\# \{ 0\leq m\leq 2d-1\, :\, \gcd (B_0+q_1 m,d)=1\}
=2\varphi(d_1)d_2 .
\end{equation*}
It follows in particular that there are always at least two
numbers $B$ as in \eqref{4.5} for which $\gcd (B,d)=1$, and as
soon as $d_2\geq 2$ or $d_1\geq 3$, there are at least four such
numbers. Since at most two numbers $B$ as in \eqref{4.5} are
divisible by $p$, and for small values of $d$ we know that at most
one number $B$ as in \eqref{4.5} is divisible by $p$, we conclude
that in all cases we have
\begin{equation}\label{4.9}
\# \{ 1\leq B\leq 2d\, :\, \gcd (B,dp)=1,\
dpb=Ba\hspace{-8pt}\pmod{q}\} \geq \varphi(d_1)d_2 \geq
\varphi(d).
\end{equation}

Combining \eqref{4.9} with \eqref{4.2} we infer that
\begin{equation*}
\# \NN_{a,b,q} \gg \frac{q_1 \varphi(d)}{\ln q_1}
=\frac{q\varphi(d)}{d\ln q_1}\geq \frac{q}{\ln q}\cdot
\frac{\varphi(d)}{d} \ ,
\end{equation*}
and the lemma is completed using the inequalities
\begin{equation*}
\frac{\varphi(d)}{d}\geq \frac{\varphi(q)}{q} \gg \frac{1}{\ln \ln
q}\ .
\end{equation*}
\qed

\medskip

Let now $q$ be a large positive integer, let $a,b\in \{
1,\dots,q\}$ such that $\gcd (a,b,q)=1$, and let $Q$ be a positive
integer larger than $q$. In our applications $Q$ will be at least
of the order of magnitude of $q^{\frac{4}{3}}$. We will construct
some sets of lattice points inside the square $[0,Q]^2$, indexed
by the set $\NN_{a,b,q}$. Precisely, we select a positive integer
$M$, which will be chosen later to be the integer part of a
certain fractional power of $Q$, and for each pair $(A,B)\in
\NN_{a,b,q}$ we construct a set $\MM_{A,B}=\MM_{A,B}(a,b,q,Q,M)$
as follows. Fix $(A,B)\in \NN_{a,b,q}$. To make a choice assume
that $B\leq A$. Let $C$ be the integer defined by the equation
\begin{equation}\label{4.10}
bA-aB=qC.
\end{equation}
Denote by $u$ the unique integer satisfying
\begin{equation*}
u=-\overline{B} \, C\pmod{A},\qquad 0\leq u\leq A-1,
\end{equation*}
where $\overline{B}$ is the multiplicative inverse of $B$ modulo
$A$, and put $v=\frac{Bu+C}{A}\, .$ Then $v$ is an integer. Also,
from the inequalities
\begin{equation*}
-B\leq -\frac{aB}{q}<\frac{bA-aB}{q}=C<\frac{bA}{q}\leq A
\end{equation*}
it follows that $\frac{C}{A} \in \big( -\frac{B}{A},1\big)
\subseteq (-1,1)$. Hence
\begin{equation*}
-1<\frac{C}{A}\leq \frac{Bu+C}{A}\leq \frac{Bu}{A}+1<B+1,
\end{equation*}
so that $v\in (-1,B+1)$. Thus $0\leq v\leq B$, since $v$ is an
integer. Let now $s=\big[ \frac{Q}{A}\big]$, and define
$\MM_{A,B}$ to be the set of lattice points given by
\begin{equation*}
\MM_{A,B}=\{ (u+mA,v+mB)\, :\, s-M\leq m\leq s-1\}.
\end{equation*}
Note that the case $A=B$ can only occur when $a=b$, and in this
situation we also get $C=0$, $u=v=0$, and so $\MM_{0,0}=\{
(mA,mA)\, :\, s-M\leq m\leq s-1\}$. We have constructed $\#
\NN_{a,b,q}$ sets of the form $\MM_{A,B}$, each set $\MM_{A,B}$
consisting of $M$ lattice points. Note that $u,v,s$ in the
definition of $\MM_{A,B}$ depend on the pair $(A,B)$. In what
follows we assume that $M$ satisfies the inequality
\begin{equation}\label{4.11}
M\leq \bigg[ \frac{Q}{4q}\bigg].
\end{equation}

Define also
\begin{equation*}
\MMM_{a,b,q}=\bigcup\limits_{(A,B)\in \NN_{a,b,q}} \hspace{-15pt}
\MM_{A,B} .
\end{equation*}
Some properties of these sets are collected in the following
lemma.

\medskip

\begin{lem}\label{L4.3}
\mbox{\rm (i)} $\dist \big([0,1]^2,\MMM_{a,b,q}\big)\geq
\displaystyle \frac{Q}{3}\, .$

\mbox{\rm (ii)} $\MMM_{a,b,q} \subseteq [0,Q]^2$.

\mbox{\rm (iii)} The sets $\MM_{A,B}$ are disjoint.
\end{lem}

{\sl Proof.} (i) Owing to \eqref{4.11} and to the inequality
\begin{equation*}
[x]-\bigg[ \frac{x}{2}\bigg] \geq \frac{x}{2}-1,
\end{equation*}
we have for any $(A,B)\in \NN_{a,b,q}$ and any point
$(u+mA,u+mB)\in \MM_{A,B}$
\begin{equation*}
u+mA\geq mA\geq (s-M)A\geq \bigg( \bigg[ \frac{Q}{A}\bigg]-\bigg[
\frac{Q}{4q}\bigg] \bigg) A\geq \bigg( \bigg[ \frac{Q}{A}\bigg]
-\bigg[ \frac{Q}{2A}\bigg] \bigg) A\geq \frac{Q}{2}-A\geq
\frac{Q}{2}-2q.
\end{equation*}
Recall that $Q$ is much larger than $q$. It follows that the
distance between any two points $P\in [0,1]^2$ and $P^\prime \in
\MMM_{a,b,q}$ satisfies
\begin{equation}\label{4.12}
\| PP^\prime \| \geq \frac{Q}{2}-2q-1\geq \frac{Q}{3}\ .
\end{equation}

(ii) For any $(A,B)\in \NN_{a,b,q}$ with, say, $A\geq B$, and any
point $P=(u+mA,v+mB)\in \MM_{A,B}$, one has $u+mA\leq u+(s-1)A\leq
sA\leq Q$. Also, $0\leq v+mB\leq v+(s-1)B\leq sB=\big[
\frac{Q}{A}\big] B\leq Q$, since $B\leq A$. Hence all points $P\in
\MMM_{a,b,q}$ lie inside the square $[0,Q]^2$.

(iii) Assume that there is a lattice point $P=(n_1,n_2)$ which
belongs to two sets $\MM_{A,B}$ and $\MM_{A^\prime,B^\prime}$ with
$(A,B),(A^\prime,B^\prime)\in \NN_{a,b,q}$ and $(A,B)\neq
(A^\prime,B^\prime)$. Assume first that $B\leq A$ and $B^\prime
\leq A^\prime$. Then
\begin{equation}\label{4.13}
n_1=u+mA,\qquad n_2=v+mB
\end{equation}
for some $m\in \{ s-M,\dots,s-1\}$, and similarly
\begin{equation*}
n_1=u^\prime+m^\prime A^\prime,\qquad n_2=v^\prime+m^\prime
B^\prime
\end{equation*}
for some $m^\prime \in \{ s^\prime -M,\dots,s^\prime -1\}$, with
$u,v,s,u^\prime,v^\prime,s^\prime$ given by appropriate
definitions.

We compute the ratio $\frac{qn_2-b}{qn_1-a}$ in two ways. First,
by \eqref{4.13}, \eqref{4.10} and the equality $v=\frac{Bu+C}{A}$
we have
\begin{equation*}
\frac{qn_2-b}{qn_1-a}=\frac{qv+qmB-b}{qu+qmA-a}
=\frac{qBu+qC+AqmB-Ab}{A(qu+qmA-a)}
=\frac{qBu-aB+AqmB}{A(qu+qmA-a)}=\frac{B}{A}\ .
\end{equation*}
By a similar computation we also have
\begin{equation*}
\frac{qn_2-b}{qn_1-a}=\frac{B^\prime}{A^\prime} =\frac{B}{A}\ .
\end{equation*}
Since $A,B,A^\prime,B^\prime$ are all positive and $\gcd
(B,A)=\gcd (B^\prime,A^\prime)=1$, this forces $A=A^\prime$ and
$B=B^\prime$.

In general we get by the same argument $\min \{ A,B\}=\min \{
A^\prime,B^\prime\}$ and $\max \{ A,B\}=\max
\{A^\prime,B^\prime\}$, thus $(A^\prime,B^\prime)\in \{
(A,B),(B,A)\}$. If $A^\prime=B$ and $B^\prime=A$, we get
$\frac{A}{B}=\frac{B}{A}$, hence $A^2=B^2$ and $A=B=1$. \qed

\bigskip

Fix now a point $P_{(x,y)}\in [0,1]^2$ and a pair $(A,B)\in
\NN_{a,b,q}$. Also, choose any two points $P,P^\prime \in
\MM_{A,B}$, say $P=(u+mA,v+mB)$ and $P^\prime =(u+m^\prime
A,v+m^\prime B)$, with $m,m^\prime\in \{ s-M,\dots,s-1\}$.
Consider the angle $\theta=\angle P^\prime P_{(x,y)}P$. Then by
\eqref{4.12}
\begin{equation}\label{4.14}
\begin{split}
\vert \sin \theta \vert & =\frac{2\area \triangle P^\prime
P_{(x,y)}P}{\| PP_{(x,y)}\| \, \| P^\prime P_{(x,y)}\|} \leq
\frac{18 \area \triangle P^\prime P_{(x,y)}P}{Q^2} \\
& =\frac{9\vert (u+m^\prime A-x)(v+mB-y)-(u+mA-x)(v+m^\prime
B-y)\vert}{Q^2} \\
& =\frac{9\vert m^\prime -m\vert \cdot \vert
A(v-y)-B(u-x)\vert}{Q^2} \\
& \leq \frac{9M\vert A(v-y)-B(u-x)\vert}{Q^2} \ .
\end{split}
\end{equation}
Next, using the equality $v=\frac{Bu+C}{A}$, we rewrite
\eqref{4.14} as
\begin{equation}\label{4.15}
\vert \sin \theta \vert \leq \frac{9M\vert C+Bx-Ay\vert}{Q^2} \ .
\end{equation}
Since $\vert C\vert \leq \max \{ A,B\} \leq 2q$ and $0\leq x,y\leq
1$, we see that $\vert C+Bx-Ay\vert \ll q$. If $M$ satisfies
\eqref{4.11}, then $M\vert C+Bx-Ay\vert \ll Q$. As a consequence
of \eqref{4.15} we also have
\begin{equation*}
\vert \sin \theta \vert \ll \frac{1}{Q} \ .
\end{equation*}
From \eqref{4.15} we also infer that
\begin{equation}\label{4.16}
\vert \theta \vert \leq \frac{\pi\vert \sin \theta\vert}{2} \ll
\frac{M\vert C+Bx-Ay\vert}{Q^2}\ .
\end{equation}

By \eqref{4.10} and \eqref{4.16} we derive that
\begin{equation}\label{4.17}
\vert \theta \vert \ll \frac{M\vert Cq+Bqx-Aqy\vert}{qQ^2}
=\frac{M\vert A(b-qy)+B(qx-a)\vert}{qQ^2} \ .
\end{equation}
As a consequence of \eqref{4.17} and of the inequalities $1\leq
A,B\leq 2q$, we have
\begin{equation}\label{4.18}
\vert \theta\vert \ll \frac{M(\vert b-qy\vert+\vert
qx-a\vert)}{Q^2} \ ,
\end{equation}
uniformly for all pairs $(A,B)\in \NN_{a,b,q}$ and all pairs of
points $P,P^\prime \in \MM_{A,B}$.

\medskip

{\sl Proof of Theorem \ref{T1.2}.} Fix $(x,y)\in [0,1]^2$,
$\lambda=(\lambda_1,\dots,\lambda_5)\in \R^5$ with
$\lambda_1,\dots,\lambda_5>0$, and $\delta >0$. The case when both
$x$ and $y$ are rational numbers is clear. In this case, if we fix
an integer $m_0\geq 1$ for which both $m_0x$ and $m_0y$ are
integers, and consider the set of lattice points $\mathcal{A}=\{
(\ell m_0x,\ell m_0 y)\, :\, \ell=1,2,\dots,\big[
\frac{Q}{m_0}\big]\big\}$, then all these points lie on the same
line, that passes through $P_{(x,y)}$. Then all the $6$-tuples of
distinct elements from $\mathcal{A}$ will contribute to
$\RR^{(6)}_{(x,y),Q}(\lambda)$. Since $\# \mathcal{A}=\big[
\frac{Q}{m_0}\big]$, it follows that
\begin{equation*}
\RR^{(6)}_{(x,y),Q}(\lambda) \gg \frac{1}{\# \Box_Q}\cdot
\frac{Q^6}{m_0^6} \gg \frac{Q^4}{m_0^6} \ .
\end{equation*}
Note that in this case the $3$-level correlations already diverge
as $Q\rightarrow \infty$.

Consider now the case when at least one of $x,y$ is irrational.
With $x,y,\lambda$ and $\delta$ fixed, choose a large positive
integer $Q$. Let $1<T<Q$ be a parameter, whose precise value will
be chosen later and will be the integer part of a fractional power
of $Q$. By Minkowski's convex body theorem (see
\cite[Thm.\,6.25]{NZM} for the formulation used here), there
exists an integer $1\leq q\leq T$ for which
\begin{equation}\label{4.19}
\langle qx\rangle \leq \frac{1}{\sqrt{T}}\qquad \mbox{\rm and}
\qquad \langle qy\rangle \leq \frac{1}{\sqrt{T}} \ ,
\end{equation}
where $\langle \, \cdot \, \rangle$ denotes here the distance to
the closest integer. Let $a$ and $b$ denote the closest integers
to $qx$ and $qy$ respectively. Then $0\leq a,b\leq q$, $\max \{
a,b\}>0$, and \eqref{4.19} gives
\begin{equation}\label{4.20}
\vert qx-a\vert \leq \frac{1}{\sqrt{T}} \qquad \mbox{\rm
and}\qquad \vert qy-b\vert \leq \frac{1}{\sqrt{T}} \ .
\end{equation}
Dividing if necessary $a$, $b$ and $q$ by $\gcd(a,b,q)$, we may
assume in what follows that $\gcd(a,b,q)=1$.

We will have $T\rightarrow \infty$ as $Q\rightarrow \infty$, and
since at least one of $x,y$ is irrational, this forces
$q\rightarrow \infty$ as $Q\rightarrow \infty$. Then all our
previous results valid for large $q$ are applicable.

Let $M$ be a positive integer satisfying \eqref{4.11}, whose
precise order of magnitude will be chosen later. Consider the
disjoint subsets $\MM_{A,B}$ of $\Box_Q$, with $(A,B)\in
\NN_{a,b,q}$. By \eqref{4.18} we know that for any $(A,B)\in
\NN_{a,b,q}$ and any $P,P^\prime\in \MM_{A,B}$, the measure of the
angle $\angle PP_{(x,y)}P^\prime$ satisfies
\begin{equation}\label{4.21}
\vert \theta_{P,P^\prime}\vert \ll \frac{M(\vert b-qy\vert+\vert
qx-a\vert)}{Q^2} \ .
\end{equation}
Plugging \eqref{4.20} in \eqref{4.21} we find that
\begin{equation*}
\vert \theta_{P,P^\prime}\vert \ll \frac{M}{Q^2 \sqrt{T}} \ .
\end{equation*}

If we take the order of magnitude of $M$ to be slightly smaller
than that of both $\frac{Q}{q}$ and $\sqrt{T}$, for instance
\begin{equation}\label{4.22}
M=\min \bigg\{ \bigg[ \frac{Q}{4q}\bigg],\bigg[
\frac{\sqrt{T}}{\ln Q}\bigg]\bigg\} ,
\end{equation}
then $M$ will satisfy \eqref{4.11} on the one hand, and on the
other hand we will have
\begin{equation*}
\vert \theta_{P,P^\prime} \vert \ll \frac{1}{Q^2 \ln Q}\ .
\end{equation*}
It follows that for $Q$ large enough in terms of
$\lambda_1,\dots,\lambda_5$, all the $6$-tuples $(P_1,\dots,P_6)$
of distinct points from $\MM_{A,B}$ will contribute to
$\RR^{(6)}_{(x,y),Q}(\lambda)$. Therefore, since $\# \MM_{A,B} =M$
for each $(A,B)\in \NN_{a,b,q}$, we derive that
\begin{equation*}
\RR^{(6)}_{(x,y),Q}(\lambda) \geq \frac{1}{N}\, \cdot \#
\NN_{a,b,q} \cdot M(M-1)\cdots (M-5).
\end{equation*}
Here $N=\# \Box_Q =(2Q+1)^2$, and from \eqref{4.22} and Lemma
\ref{L4.2} it follows that
\begin{equation}\label{4.23}
\RR^{(6)}_{(x,y),Q}(\lambda) \gg \frac{M^6 q}{Q^2 \ln q\ln \ln q}
\gg \min \bigg\{ \frac{Q^4}{q^5 \ln^2 q} \, ,\ \frac{T^3 q}{Q^2
\ln^8 Q} \bigg\} .
\end{equation}

We now choose $T=[Q^{\frac{3}{4}}]$. Then, no matter how small $q$
might be, we have
\begin{equation}\label{4.24}
\frac{T^3 q}{Q^2 \ln^8 Q} \geq \frac{Q^{\frac{1}{4}}}{\ln^8 Q}
>Q^{\frac{1}{4}-\delta} \qquad \mbox{\rm for large $Q$.}
\end{equation}
Also, since $q\leq T$, it follows that
\begin{equation}\label{4.25}
\frac{Q^4}{q^5 \ln^2 q} \geq \frac{Q^4}{T^5 \ln^2 T} \geq
\frac{Q^{\frac{1}{4}}}{\ln^2 Q} >Q^{\frac{1}{4}-\delta} \qquad
\mbox{\rm for large $Q$.}
\end{equation}

Now \eqref{IV.28} follows from \eqref{4.23}, \eqref{4.24} and
\eqref{4.25}. \qed

\bigskip

\section{Appendix}

For a fixed integer $q\geq 2$, we consider for any integers
$m,n,h$ and any sets $I,I_1,I_2 \subset \R$ the Kloosterman sum
\begin{equation*}
K(m,n;q)=\sum\limits_{\substack{x\hspace{-8pt}\pmod{q} \\ \gcd
(x,q)=1}} e\bigg( \frac{mx+n\bar{x}}{q}\bigg),
\end{equation*}
the incomplete Kloosterman sum
\begin{equation*}
K(m,n;q)=\sum\limits_{\substack{x\in I \\ \gcd (x,q)=1}} e\bigg(
\frac{mx+n\bar{x}}{q}\bigg),
\end{equation*}
and the set
\begin{equation*}
N_{q,h} (I_1,I_2)=\big\{ (x,y)\in I_1 \times I_2\, :\, \gcd
(x,q)=1,\ xy=h\hspace{-7pt}\pmod{q}\big\} .
\end{equation*}
Here $\bar{x}$ denotes the multiplicative inverse of $x\pmod{q}$.

\medskip

\begin{lem}\label{LA1} {\em (\cite[Lemma 2.2]{BCZ})}
Suppose that $f$ is a piecewise $C^1$ function on $[a,b]$. Then
\begin{equation*}
\sum\limits_{\substack{a<k\leq b \\ \gcd (k,q)=1}} f(k)
=\frac{\varphi(q)}{q} \int_a^b f+O\Big( q^\eps ( \|
f\|_\infty+V_a^b f) \Big).
\end{equation*}
\end{lem}

\begin{lem}\label{LA2} {\em (\cite[Lemma 1.6]{BCZ})}
For any interval $I\subset [0,q)$ and any integer $n$
\begin{equation*}
\vert S_I (0,n,q)\vert \ll \gcd (n,q)^{\frac{1}{2}}
q^{\frac{1}{2}+\eps} .
\end{equation*}
\end{lem}

{\sl Proof.} We write
\begin{equation}\label{A1}
\begin{split}
S_I(0,n,q) & =\sum\limits_{\substack{x\hspace{-8pt} \pmod{q} \\
\gcd (x,q)=1}} e\bigg( \frac{n\bar{x}}{q}\bigg) \sum\limits_{y\in
I} \frac{1}{q} \sum\limits_{k=1}^q e\bigg( \frac{k(y-x)}{q}\bigg)
\\ & =\frac{1}{q} \sum\limits_{k=1}^q \sum\limits_{y\in I}
e\bigg( \frac{ky}{q}\bigg)
\sum\limits_{\substack{x\hspace{-8pt} \pmod{q}
\\ \gcd (x,q)=1}} e\bigg( \frac{-kx+n\bar{x}}{q} \bigg)  \\
& =\frac{1}{q} \sum\limits_{k=1}^q \sum\limits_{y\in I} e\bigg(
\frac{ky}{x}\bigg) K(-k,n,q).
\end{split}
\end{equation}
The inner sum is a geometric progression and can be bounded as
\begin{equation}\label{A2}
\Bigg| \sum\limits_{y\in I} e\bigg( \frac{ky}{q}\bigg) \Bigg| \leq
\min \left\{ \vert I\vert ,\frac{1}{2\big\| \frac{k}{q}\big\|}
\right\} .
\end{equation}
By \eqref{A1} and \eqref{A2}
\begin{equation*}
\begin{split}
\vert S_I (0,n,q)\vert & \leq \frac{1}{q}\ K(0,n,q) \vert I\vert
+\frac{1}{q} \Bigg| \sum\limits_{k=1}^{q-1} \frac{1}{2\big\|
\frac{k}{q}\big\|} \ K(-k,n,q)\Bigg| \\
& \ll \frac{q q^{\frac{1}{2}+\eps} \gcd (n,q)^{\frac{1}{2}}}{q}
+\frac{q^{\frac{1}{2}+\eps} \gcd (n,q) q\log q}{q} \ll \gcd
(n,q)^{\frac{1}{2}} q^{\frac{1}{2}+2\eps} .
\end{split}
\end{equation*}
\qed

\medskip

\begin{prop}\label{PA1}
Suppose that $I_1,I_2 \subset [0,q)$ are intervals. Then for any
integer $h$
\begin{equation*}
N_{q,h}(I_1,I_2)=\frac{\varphi (q)}{q^2} \ \vert I_1\vert\, \vert
I_2 \vert+O_\delta \big( q^{\frac{1}{2}+\delta} \gcd (h,q)\big).
\end{equation*}
\end{prop}

{\sl Proof.} We write
\begin{equation*}
N_{q,h}(I_1,I_2)=\frac{1}{q} \sum\limits_{\substack{x\in I_1,y\in
I_2 \\ \gcd (x,q)=1}} \sum\limits_{k=0}^{q-1} e\bigg(
\frac{k(y-h\bar{x})}{q}\bigg) =M+E,
\end{equation*}
where the main term can be expressed, as a result of Lemma
\ref{LA1}, as
\begin{equation*}
M =\frac{1}{q} \sum\limits_{\substack{x\in I_1,y\in I_2 \\ \gcd
(x,q)=1}} 1 =\frac{1}{q} \bigg( \frac{\varphi(q)}{q}\ \vert I_1
\vert+O(q^\eps)\bigg) \big( \vert I_2 \vert +O(1)\big)
=\frac{\varphi (q)}{q^2}\ \vert I_1 \vert \, \vert I_2
\vert+O(q^\eps).
\end{equation*}
The error is given by
\begin{equation*}
E =\frac{1}{q} \sum\limits_{k=1}^{q-1} \sum\limits_{y\in I_2}
e\bigg( \frac{ky}{q}\bigg) \sum\limits_{\substack{x\in I_1 \\ \gcd
(x,q)=1}} e\bigg( -\frac{hk\bar{x}}{q}\bigg) =\frac{1}{q}
\sum\limits_{k=1}^{q-1} \sum\limits_{y\in I_2} e\bigg(
\frac{ky}{x}\bigg) S_{I_1} (0,-hk,q).
\end{equation*}
Owing to \eqref{A2}, Lemma \ref{LA2}, and to the inequality $\gcd
(kh,q)\leq \gcd (h,q)\gcd (k,q)$, we have
\begin{equation*}
\vert E\vert \leq \frac{\gcd (h,q)^{\frac{1}{2}}}{q}\
\sum\limits_{k=1}^{q-1} \frac{1}{2\big\| \frac{k}{q}\big\|} \ \gcd
(k,q)^{\frac{1}{2}} q^{\frac{1}{2}+\eps}\leq \gcd
(h,q)^{\frac{1}{2}} q^{\frac{1}{2}+\eps}
\sum\limits_{k=1}^{\frac{q-1}{2}} \frac{1}{k}\ \gcd
(k,q)^{\frac{1}{2}} .
\end{equation*}
Writing $k=dm$, with $d=\gcd (k,m)$, we eventually get
\begin{equation*}
\vert E\vert \leq \gcd (h,q)^{\frac{1}{2}} q^{\frac{1}{2}+\eps}
\sum\limits_{d\vert q} \frac{1}{d^{\frac{1}{2}}}
\sum\limits_{m=1}^{\big[ \frac{q}{d}\big]} \frac{1}{m} \ll \gcd
(h,q)^{\frac{1}{2}} q^{\frac{1}{2}+2\eps} \log q\ll \gcd
(h,q)^{\frac{1}{2}} q^{\frac{1}{2}+3\eps}.
\end{equation*}
\qed

\begin{prop}\label{PA2}
Assume that $q\geq 1$ and $h$ are two given integers, $\II$ and
$\JJ$ are intervals of length lesser than $q$, and $f:\II \times
\JJ \rightarrow \R$ is a $C^1$ function. Then for any integer
$T>1$ and any $\delta >0$
\begin{equation*}
\begin{split}
\sum_{\substack{a\in \II,\, b\in \JJ \\
ab=h\hspace{-8pt}\pmod{q}\\ \gcd (b,q)=1}} \hspace{-3pt} & f(a,b)
=\frac{\varphi (q)}{q^2} \iint\limits_{\II \times \JJ} f(x,y) dx
dy  \\ & +O_\delta \bigg( T^2 \| f\|_\infty q^{\frac{1}{2}+\delta}
\gcd(h,q) +T\| Df\|_\infty q^{\frac{3}{2}+\delta} \gcd (h,q)+
\frac{\vert \II \vert \, \vert \JJ \vert \, \| Df\|_\infty}{T}
\bigg).
\end{split}
\end{equation*}
\end{prop}

The proof is identical to that of Lemma 2.2 in \cite{BGZ}, and
relies on Proposition \ref{PA1}.

\end{document}